\UseRawInputEncoding
\documentclass[12pt, reqno]{amsart}
\usepackage[margin=1in]{geometry}
\usepackage[mathscr]{eucal}
\usepackage{bm}
\usepackage{mathptmx}
\usepackage{amssymb}
\usepackage{amsthm}
\usepackage{dcolumn}
\usepackage[all]{xy}
\usepackage{enumitem}
\usepackage{xcolor}
\usepackage[hidelinks]{hyperref}
\usepackage{mathtools} 
\usepackage{bm} 

\def\textmatrix#1&#2\\#3&#4\\{\bigl({#1 \atop #3}\ {#2 \atop #4}\bigr)}
\def\dispmatrix#1&#2\\#3&#4\\{\left({#1 \atop #3}\ {#2 \atop #4}\right)}

\newcommand{\beg}{\begin{equation}}
        \newcommand{\eeg}{\end{equation}}
\newcommand{\ben}{\begin{eqnarray*}}
        \newcommand{\een}{\end{eqnarray*}}

\newtheorem{thm}{Theorem}[section]
\newtheorem{cor}[thm]{Corollary}
\newtheorem{lem}[thm]{Lemma}

\newtheorem{prop}[thm]{Proposition}
\numberwithin{equation}{section} \theoremstyle{definition}
\newtheorem{defn}[thm]{Definition}

\newtheorem{note}[thm]{Note}

\def\textmatrix#1&#2\\#3&#4\\{\bigl({#1 \atop #3}\ {#2 \atop #4}\bigr)}
\def\dispmatrix#1&#2\\#3&#4\\{\left({#1 \atop #3}\ {#2 \atop #4}\right)}

\newcommand{\C}{\mathbb{C}}
\newcommand{\D}{\mathbb{D}}
\newcommand{\N}{\mathbb{N}}

\newcommand{\X}{\mathbb{X}}
\newcommand{\Y}{\mathbb{Y}}
\newcommand{\Z}{\mathbb{Z}}
\newcommand{\W}{\mathbb{W}}

\newcommand{\SR}{\mathsf{R}}
\newcommand{\HS}{\mathcal{H}}

\title[Decomposition, BCL model and unitary extension of Banach space Wold isometries]{On S\l oci\'nski-type decomposition, BCL model and unitary extension of commuting Banach space Wold isometries}

\author[Jana and Pal]{Swapan Jana and Sourav Pal}

\address[Swapan Jana]{Mathematics Department, Indian Institute of Technology Bombay, Powai, Mumbai - 400076, India.} \email{ swapan.jana@iitb.ac.in , swapan.math2015@gmail.com}

\address[Sourav Pal]{Mathematics Department, Indian Institute of Technology Bombay, Powai, Mumbai - 400076, India.} \email{souravpal@iitb.ac.in , souravmaths@gmail.com}
	
\author[Saikat Roy]{Saikat Roy}
\address{Department of Mathematics, School of Advanced Science VIT-AP University, Amaravathi 522241, India} \email{saikat.roy@vitap.ac.in , saikatroy.cu@gmail.com}

	

\thanks{The first named author is supported by the ``Prime Minister's Research Fellowship (PMRF)" with Award No. PMRF-1302045. The second named author is supported in part by Core Research Grant with Award No. CRG/2023/005223 from Anusandhan National Research Foundation (ANRF) of Govt. of India. The third named author was supported by the Institute Post-doctoral Fellowship of IIT Bombay during the course of the work.}

\begin{document}

\begin{abstract}

The adjoint of an operator
plays the central role in Hilbert space operator theory. In particular, dilations
of Hilbert space contractions, extension of commuting Hilbert space isometries to
Hilbert space unitaries, orthogonal decompositions of contractions and functional
models for commuting isometries are heavily dependent on the
presence of adjoint. In this paper, we develop a theory of Banach space Wold isometries with right-complemented ranges by showing explicit
constructions of Berger-Coburn-Lebow model and unitary extension for such isometric tuple without the help of adjoint. We also show that a S\l oci\'nski type decomposition is possible for such a commuting tuple of Banach space isometries. Each of our results generalizes the corresponding result in the Hilbert space theory.
\end{abstract}

\maketitle

\section{Introduction}\label{sec:001}

\noindent
We begin with a few basic notations and terminologies that we shall follow throughout the paper. All Banach and Hilbert spaces are over the field of complex numbers $\C$. The set of all bounded operators on a Banach space $\X$ is denoted by $\mathcal B(\X)$. The \textit{identity operator} on a Banach space $\X$ will be denoted by $I_{\X}$ or often just by $I$. For two Banach spaces $\X$ and $\Y$, an \textit{isometry} is a linear map $V:\, \X \rightarrow \Y$ that satisfies $\|Vx\|=\|x\|$ for all $x\in \X$. A \textit{Banach space unitary} is a surjective isometry. If $\X_1,\X_2$ are any two normed spaces, then $\X_1 \oplus \X_2$ denotes the vector space direct sum of $\X_1, \X_2$ that consists of vectors of the form $(x_1,x_2)$ or $x_1 \oplus_2 x_2$, where $x_1 \in \X_1, \, x_2 \in \X_2$. Their \textit{direct $2$-sum} (or, \textit{orthogonal $2$-sum}) $\X_1 \oplus_2 \X_2$ is the normed linear space consisting of the pairs $(x_1,x_2)$ with $x_i\in \X_i$ ($i=1,2$) and is equipped with the norm $\|(x_1,x_2)\| = (\|x_1\|^2 + \|x_2\|^2)^{\frac{1}{2}}$. To maintain uniformity, the direct sum of two Hilbert spaces $\HS_1$ and $\HS_2$ (which is also orthogonal $2$-sum in case of Hilbert spaces) will also be denoted by $\HS_1 \oplus_2 \HS_2$. For any Banach space $\X$, the space $\ell_2(\X)$ consists of all square-summable sequences $\{x_n\}$ with entries in $\X$ and is equipped with the norm $\|\{x_n\}\|=(\sum_{n=1}^\infty \|x_n\|^2)^{\frac{1}{2}}$. For a Hilbert space contraction $T$, the defect operator and the defect space of $T$ are $D_T= (I -T^*T)^{\frac{1}{2}}$ and $\mathfrak{D}_T=\overline{Ran } \ D_T$, respectively.

\smallskip

\subsection{S\l oci\'nski-type decomposition for Banach space Wold isometries} An isometry $V$ on a Hilbert space $\HS$ orthogonally splits into two parts of which one is a unitary and the other is a unilateral shift. This is a famous result due to Wold and von Neumann which is stated below.

\begin{thm}[\cite{Wold}, von Neumann-Wold]
Let V be an isometry on a Hilbert space $\HS$. Then $\HS$ decomposes into an orthogonal sum $\HS = \HS_0 \oplus_2 \HS_1$ such that $\HS_0$ and $\HS_1$ reduce $V$ and that $V|_{\HS_0}$
is a unitary and $V|_{\HS_1}$ is a pure isometry, i.e., a unilateral shift. This decomposition is uniquely determined; indeed we have
\[
\HS_0 = \bigcap_{n=0}^\infty V^n\HS \text{ and } \HS_1 = M_+(\mathcal{L}), \quad \text{ where } \quad  \mathcal{L} = \HS \ominus V\HS \text{ and } M_+(\mathcal{L}) = \bigoplus_{n=0}^\infty V^n\mathcal{L}.
\]
\end{thm}

Note that a \textit{pure isometry} ${V}_1$ on a Hilbert space $\widehat{\HS}$ is an isometry that has no unitary part, i.e. there is no nonzero closed subspace of $\widehat{\HS}$ which is reduced by ${V}_1$ and ${V}_1$ restricted to that subspace is a unitary. A pure isometry on a Hilbert space is unitary equivalent to a unilateral shift. 

\smallskip

M. S\l oci\'nski \cite{Slocinski} found an appealing generalization of Wold decomposition of a Hilbert space isometry. In general, a pair of commuting isometries on a Hilbert space does not admit a Wold-type decomposition as was shown in the same paper, see \cite[Example 1]{Slocinski}. For this reason, S\l oci\'nski considered a pair of doubly commuting isometries on a Hilbert space and established a Wold-type decomposition for such a pair. We state the result below.

\begin{thm}[S\l oci\'nski, \cite{Slocinski}]\label{thm:102}
Suppose $V=(V_1, V_2)$ be two doubly commuting isometries on a Hilbert space $\HS.$ Then there is unique decomposition
\[
   \HS = \HS_{uu} \oplus_2 \HS_{us} \oplus_2 \HS_{su} \oplus_2 \HS_{ss},
\]
where $\HS_{ij}$ is reducing subspace of $V_1, V_2$ for all $i,j= s,u.$ Moreover, $V_1$ on $\HS_{ij}$ is unitary if $i=u,$ and is a unilateral shift if $i=s$. Similarly, $V_2$ on $\HS_{ij}$ is unitary if $j=u$ and is a unilateral shift if $j=s$.
\end{thm}

Several independent proofs of Theorem \ref{thm:102} can be found in the literature, e.g., \cite{CPS, KM, KO}. S\l oci\'nski’s result was further generalized for any number of doubly commuting isometries with finite dimensional wandering subspaces in \cite{BKS}. Also, the recent articles \cite{Pal-Sourav-decomp}, \cite{P-S-Tm} discuss in detail various decompositions of Hilbert space contractions in commutative and non-commutative settings, respectively.

\smallskip

The main hindrance in studying such results in Banach space setting is dealing with orthogonal complement of a closed subspace. The idea of orthogonal complement of a closed subspace in a Hilbert space was generalized to Banach spaces through the notion of $1$-complemented subspace. A closed subspace of a Banach space is said to be \textit{orthogonally complemented} or $1$-\textit{complemented} or \textit{right-complemented} \cite{JPR, R 2} if it is the range of a norm-one projection. However, unlike every closed subspace of a Hilbert space, a closed subspace of a Banach space may not be $1$-complemented. In fact, the Banach space isometries that have $1$-complemented ranges are not fully known till date even though there is an extensive research on this topic, e.g. see \cite{CFS, JPR} and the relevant references therein.

\smallskip

In \cite{Faulkner Huneycutt}, Faulkner and Huneycutt first introduced the notion of unilateral shift (see Definition \ref{defn:Banach-Unilateral}) on a Banach space and obtained a Wold-type decomposition for an isometry having $1$-complemented range and acting on a reflexive and smooth Banach space. This result was further generalized to reflexive Banach spaces removing the smoothness assumption by Campbell, Faulkner and Sine in \cite[Corollary 1]{CFS}. For a further reading in this direction, a reader is referred to \cite{JPR}, and the references therein. By a Wold isometry on a Banach space one naturally means an isometry that admits a Wold-type decomposition and consequently Wold isometry in Banach space setting was defined in the following way, see \cite{CFS, Faulkner Huneycutt, R 2} for details.

\begin{defn}\label{Definition:Wold-isometry-003}
Let $V$ be an isometry on a Banach space $\mathbb{X}$. Then $V$ is said to be a \textit{Wold isometry}, if $\X$ admits a decomposition $\X=\X_u\oplus \X_s$ into invariant subspaces $\X_u, \X_s$ of $V$ such that $V|_{\X_u}$ is a unitary (i.e. a surjective isometry) and $V|_{\X_s}$ is a pure isometry. Here $\X= \X_u \oplus \X_s$ means that every element $x\in \X$ can be uniquely expressed as $x= x_1 + x_2 $ for some $x_1 \in \X_u$ and $x_2 \in \X_s$. A Wold isometry $V$ is written as $V=V|_{\X_u} \oplus V|_{\X_s}$, which means $Vx=V|_{\X_u}x_1 \oplus V|_{\X_s}x_2 \ $, and this expression is called the \textit{Wold decomposition} of $V$.
\end{defn}

Note that for a Wold isometry $V$ on $\X$, the spaces $\X_u, \X_s$ are uniquely determined as
\begin{equation}\label{Wold decomposition Banach spaces}
\X_u= \bigcap_{n= 0}^{\infty} V^n\mathbb{X}, \quad \X_s= \bigoplus_{n=0}^{\infty}V^n\mathcal{L},
\end{equation}
where $\mathcal L$ is vector space complement of $V\X$.
Unlike pure isometries on a Hilbert space, a pure Banach space isometry may not always be a unilateral shift as was shown in \cite{PAF}. 

\smallskip

For an isometry $V$ on a Banach space $\X$, an operator $\Delta \in \mathcal B(\X)$ with $\|\Delta\|=1$ is called a \textit{left-inverse} of $V$ if $\Delta V =I_{\X}$. Note that left-inverse of a Banach space isometry $V$ may not be unique. 

\smallskip

\noindent \textbf{Our contribution.} In order to find an appropriate generalization of S\l oci\'nski's result in Banach space setting, first we need to have a meaningful definition of doubly commuting Banach space isometries. Recall that a commuting tuple of operators $(T_1, \dots, T_n)$ on a Hilbert space is called \textit{doubly commuting} if $T_i^*T_j=T_jT_i^*$ for $1\leq i,j \leq n$ and $i \neq j$. For a Hilbert space isometry $V$, the adjoint $V^*$ plays the role of left-inverse of $V$ and hence double commutativity of a tuple of commuting Hilbert space isometries can be defined in terms of left-inverses of isometries. Though Banach space operators do not have adjoints in general, but a Banach space Wold isometry with right-complemented range can have a left-inverse. Moreover, if the Banach space is reflexive, smooth and strictly convex, then left-inverse of such an isometry always exists and is unique as is shown in Lemma \ref{lem:202}. We capitalize this fact and define double commutativity for such Banach space isometries in terms of left-inverse in the following way.
 \begin{defn}\label{defn:sec-01-01}
 Let $\X$ be a reflexive, smooth and strictly convex Banach space and let $V_1, \dots, V_n$ on $\X$ be commuting Wold isometries with right-complemented ranges. Let $\Delta_i$ be the left-inverse of $V_i$ for $1\leq i \leq n$. Then the tuple $(V_1, \dots, V_n)$ is said to be \textit{doubly commuting} if $\Delta_iV_j=V_j\Delta_i$ for $1\leq i, j \leq n$ and $i \neq j$.  
 \end{defn}
A second step towards generalizing S\l oci\'nski's result to Banach space setting is to define meaningfully a notion of reducing subspace for Banach space isometries. As we are considering only Wold isometries with right-complemented ranges acting on reflexive, smooth and strictly convex Banach spaces so that we can have left-inverses of them, we define reducing subspace only for such Banach space isometries.

\begin{defn}\label{defn:ch-10-reducing-subspace}
Let $V$ be  a Wold isometry with right-complemented range on a reflexive, smooth and strictly convex Banach space. Then a closed linear subspace $M$ is said to be a \textit{reducing subspace} of $V$ if $M$ is invariant under both $V$ and its left-inverse $\Delta$. The closed subspace $M$ is said to be a \textit{joint reducing subspace} of a family of such Wold isometries $\{V_i\}_{i\in I}$, if $M$ is a reducing subspace for every $V_i$, $i\in I$.
\end{defn} 
 
 In Theorem \ref{thm:302}, we present a S\l oci\'nski-type decomposition to a pair of doubly commuting isometries with right-complemented ranges acting on a reflexive, smooth and strictly convex Banach space. Further, following the same mechanism we generalize this result to any finite isometric tuple of same kind in Theorem \ref{thm:303}, which is stated below.
 
 \begin{thm}\label{CH1:thm:303}
Let $(V_1, \dotsc, V_m)$ be a tuple of doubly commuting Wold isometries with right-complemented ranges on a reflexive, smooth and strictly convex Banach space $\X$. Let $A=\{u,s\}$ and $B=\{1,2,\dotsc, m\}$. Then there are $2^m$ many subspaces $\left\{\X_f: f\in B^A\right\}$ such that the following holds:
\begin{enumerate}
\item[(i)] $\X_f$ is a joint reducing subspace of $\{V_i\}_{_1}^m$ for each $f\in B^A$;

\item[(ii)] for $1\leq i \leq m$, and $f\in B^A$, $V_i|_{\X_f}$ is unitary if $f(i) =u$ and unilateral shift if $f(i) =s$; 

\item[(iii)] for the constant function $f(j)=u$, where $(f\in B^A)$, $\X_f$ is the maximal joint reducing subspace of $\{V_i\}_{_1} ^m$ such that $V_{i}|_{\X_f}$ is a unitary for $1\leq i \leq m$;

\item[(iv)] for the constant function $g(j)=s$, where $(g\in B^A)$, $\X_g$ is the maximal joint reducing subspace of $\{V_i\}_{_1} ^m$ such that $V_{i}|_{\X_g}$ is a unilateral shift for $1\leq i \leq m$.

\end{enumerate}
\end{thm}

\subsection{Unitary extension of a commuting tuple of Banach space Wold isometries} Every Hilbert space isometry extends to a Hilbert space unitary, e.g., see Chapter $1$ of \cite{Nagy Foias}. More precisely, the unilateral shift part of a Hilbert space isometry naturally extends to a bilateral shift, which is a Hilbert space unitary. Bilateral shifts on Banach spaces, which were introduced by Campbell, Faulkner, and Sine in \cite[Definition 2]{CFS}, are not as well-behaved as they are on Hilbert spaces. In fact, a bilateral shift on a Banach space may not be even a bounded operator as was shown in \cite{Gellar}. Consequently, extension of a Banach space unilateral shift to a bilateral shift is not straightforward.

\smallskip

A nonempty family of commuting Hilbert space isometries admits a simultanenous extension to a commuting family of Hilbert space unitaries, e.g., see Chapter $1$ of \cite{Nagy Foias}. Therefore, possibility of such an extension in Banach space setting comes as a natural question. In \cite{Douglas}, Douglas addressed this problem for an isometric representation on a Banach space of a semigroup. An \textit{isometric} (\textit{or, unitary}) \textit{representation} on a  Banach space $\X$ of a semigroup $\Sigma$ is a map $V: \Sigma \to \mathcal{B}(\X)$ such that (i) $V(\sigma)$ is isometry (or, unitary)  for all $\sigma\in \Sigma$, (ii) $V(\sigma \tau) = V(\sigma) V(\tau)$ for all $\sigma, \tau$ in $\Sigma$. Douglas \cite{Douglas} proved that for every isometric representation $V: \Sigma \to \mathcal{B}(\X)$, there is a Banach space $\Y\supseteq \X$ and a unique unitary representation $U: \Sigma \to \mathcal{B}(\Y)$ such that $U$ extends $V$.  

\smallskip

\noindent \textbf{Our contribution.} In Section \ref{Sec:12}, we explicitly construct an extension for a finite commuting tuple of Wold isometries with right-complemented ranges acting on a reflexive, smooth and strictly convex Banach space. Our construction is analogous to the same for the Hilbert space isometric tuple (see Chapter $1$ of \cite{Nagy Foias}). For a given tuple of such commuting Banach space isometries, our extension is a commuting tuple of invertible Banach space operators that are not necessarily unitaries. Then we characterize the isometric tuples for which such extensions are commuting tuple of Banach space unitaries. This characterization involves the decomposition of the underlying Banach space in terms of $2$-direct sum complement of the range of each isometry. This is Theorem \ref{thm:503} in this paper and is stated below.
\begin{thm}\label{CH1:thm:503}
Let $n\geq 2$ and let $(V_1, \dotsc, V_n)$ be a tuple of Wold isometries with right-complemented ranges acting on a reflexive, smooth and strictly convex Banach space $\X$. Then there is a Banach space $\Y\supseteq\X$ and a tuple $(U_1, \dotsc, U_n)$ of commuting invertible operators on $\Y$ such that $U_1^{m_1} \dotsc U_n^{m_n}$ is an extension of $V_1^{m_1} \dotsc V_n^{m_n}$ for all $m_i\geq 0$, $1\leq i \leq n$. Moreover, $(U_1,\dotsc, U_n)$ is a tuple of commuting isometries if and only if $\X = V_i (\X) \oplus_2 \mathcal{L}_i$ for some closed subspaces $\mathcal{L}_i$ for $1\leq i \leq n$.
\end{thm}
Douglas's results in \cite{Douglas} show existence of extension of one or more Banach space isometries to unitaries. However, in Section \ref{Sec:12} of this paper, we explicitly construct such unitary extensions in certain cases.

\smallskip

\subsection{A Berger-Coburn-Lebow model Banach space Wold isometries} Berger, Coburn and Lebow \cite{BCL} obtained a model theorem for commuting Hilbert space isometries whose product was a pure isometry , i.e., a unilateral shift $M_z$ of some multiplicity. The fact that the commutants of $M_z$ on a vector-valued Hardy space $H^2(E)$ with $E$ being a Hilbert space, are Toeplitz operators $M_{\phi}$ with $\phi \in \mathcal H^{\infty}\left( B(E) \right)$ is used in the Berger-Coburn-Lebow (BCL) model for commuting Hilbert space isometries.

\smallskip

\noindent \textbf{Our contribution.} In Section \ref{Sec:11}, we present an explicit model for a commuting tuple of Wold isometries with right-complemented ranges acting on a reflexive, smooth and strictly convex Banach space. Naturally, our model for commuting Banach space Wold isometries is a generalization of the BCL model in Hilbert space setting. To do so, first we need to have an explicit form of the commutants of a Banach space unilateral shift. The commutant of a Banach space shift of multiplicity one was found by Crownover in \cite{RMC}, though as per our knowledge the same for a unilateral shift of arbitrary multiplicity is not known till date. For this reason, we find in Theorem \ref{thm:401} the commutants of such a shift operator with arbitrary multiplicity and capitalize this to construct an explicit model for commuting Wold isometries with right-complemented ranges in Theorem \ref{thm:411}, which is stated below.
\begin{thm}\label{CH1:thm:411}
Let $(V_1, \dotsc, V_m)$ be a tuple of commuting Wold isometries with right-complemented ranges on a reflexive, smooth and strictly convex Banach space $\X$ such that $V= \prod_{i=1}^m V_i$ is a pure isometry. Then the following assertions hold:
\begin{enumerate}
\item[(i)] There is a unique closed subspace $\W$ and a tuple $(\mathbf{P_1}, \dotsc, \mathbf{P_m})$ of norm-one projections on $\W$ such that
\begin{equation*} 
\begin{split}
  (V_1, \dotsc, V_m, V) & \equiv  \left(\mathbf{T}_{M_1 + zN_1}, \dotsc, \mathbf{T}_{M_m + zN_m}, \mathbf{T}_{zI_\W} \right) \quad \text{ on } \quad \mathcal{S},\ \ \text{ where } \\
  & \quad  M_k = V_k(I_\W - \mathbf{P_k}), \ \ N_k = \prod_{i=1,~ i\neq k}^m \Delta, \ \ 1\leq k \leq m,
\end{split}
\end{equation*}
and $\mathcal{S}$ is the sequence space $($as in $(\ref{eq:202}))$ associated with  the isometry $V$, and $\Delta_i$ is left-inverse of $V_i$.

\smallskip

\item[(ii)] The left-inverse of $\mathbf{T}_{M_k + zN_k}$ $(1\leq k\leq m)$ is again a Toeplitz block matrix and is given by $\mathbf{T}_{C_k + \bar{z}D_k}$ where
\begin{align*}
& C_k = \begin{multlined}[t]
 \sum_{i=2}^{m-2} \left(M_{\sigma^k(1)}\dotsc M_{\sigma^k(i-1)}N_{\sigma^k(i)}M_{\sigma^k(i+1)}\dotsc M_{\sigma^k(m-1)}\right) \\
+ M_{\sigma^k(1)}\dotsc M_{\sigma^k (m-2)}N_{\sigma^k (m)} + N_{\sigma^k(1)} M_{\sigma^k (2)}\dotsc M_{\sigma^k (m-1)}, 
\end{multlined} \\
& D_k = \prod_{i=1}^{m-1} M_{\sigma^k(i)}, \quad \text{where} \quad \sigma(i) = i+k, \ \ \text{for} \ \ 1\leq i \leq m-1, \ \sigma(m) =1,
\end{align*}
and summation over empty set is excluded.

\smallskip

\item[(iii)] The operators $\{\Delta_k M_k, M_k \Delta_k:~ 1\leq k \leq m\}$ on $\W$ are bounded projections of norm at most 2 and have norm 1 if and only if 
\[
  \left(\bigoplus_{i=1}^{m-2} \prod_{j=1}^i V_{\sigma^k(j)} \left(\mathcal{L}_{\sigma^k(i+1)} \right) \oplus \mathcal{L}_{\sigma^k(1)} \right) \perp_B \left(\prod_{i=1}^{m-1} V_{\sigma^k(i)}\mathcal{L}_{\sigma^k(m)}\right),~~ 1\leq k\leq m,
\]
where direct sum over empty set is excluded
\end{enumerate}
\end{thm}

\section{Preliminaries}\label{sec:02}

\noindent
In this Section, we recall from the literature a few basic concepts that will be used throughout the paper. We begin with the definition of reflexive, smooth and strictly convex Banach spaces.

\medskip

\noindent \textbf{$\bullet$ Reflexive, smooth and strictly convex Banach spaces:} Let $\X$ be a Banach space and let $\pi_\X:\X\to \X^{**}$ be the canonical embedding defined by
\[
    \pi_\X (x) = \hat{x}, \qquad \hat{x}(x^*) = x^*(x), \quad x^* \in \X^*.
\]
The map $\pi_\X$ is a linear isometry. The space $\X$ is called \textit{reflexive} if $\pi_\X$ is surjective. A linear functional $f\in \X^*$ is said to be a \textit{support functional} at $x\in \X$ if $f(x)=\|f\|\|x\|$ and $\|f\|=\|x\|$. The collection of all support functionals at $x$, denoted by $J(x)$, is the following set: 
\begin{equation}\label{eqn:new-021}
 J(x)=\{f_x\in \X^*:~\|f_x\|=\|x\|,~f_x(x)=\|x\|^2\}.
\end{equation}
It follows from Hahn-Banach Theorem that the collection $J(x)$ is nonempty.  A point $x\in \X$ is said to be a \textit{smooth point} if there is a unique support functional at $x$. The space $\X$ is said to be \textit{smooth} if every non-zero vector in $\X$ is a smooth point. The Banach space $\X$ is said to be \textit{strictly convex} or \textit{rotund} or \textit{strictly normed} if $\|tx_1+(1-t)x_2\|<1$ whenever $x_1,x_2$ are distinct unit vectors in $\X$ and $0<t<1$. It is well known that if $\X^*$ is smooth then $\X$ is strictly convex and if $\X^*$ is strictly convex then $\X$ is smooth. Therefore, a reflexive Banach space $\X$ is smooth or strictly convex if and only if $\X^*$ is strictly convex or smooth respectively.

\medskip

\noindent
\textbf{$\bullet$ Birkhoff-James orthogonality:}  Let $\X$ be a Banach space and $x,\; y\in \X$. Then $x$ is said to be orthogonal to $y$ in the sense Birkhoff-James \cite{Birkhoff, James 1}, written as $x \perp_B y$, if $\|x + \lambda y\| \geq \|x\|$ for all scalars $\lambda$. Similarly, for two subspaces $\X_1$ and $\X_2$ of $\X$, we say $\X_1 \perp_B \X_2$, if $x_1 \perp_B x_2$ for all $x_i\in \X_i$, $i=1,2$. When $\X$ is a Hilbert space, the notion of Birkhoff-James orthogonality is equivalent to the inner product orthogonality. However, unlike to the inner product orthogonality, Birkhoff-James orthogonality may not be symmetric and left-additive in general. In fact, regarding the symmetry and left-additivity of Birkhoff-James orthogonality, we have the following famous result due to James.

\begin{thm}[James, \cite{James 1}]\label{thm:201}
Let $\X$ be Banach space with $dim~ \X \geq 3$. Then the following are equivalent:
\begin{enumerate}
\item[(i)] $X$ is a Hilbert space;
\smallskip
\item[(ii)] Birkhoff-James orthogonal is symmetric, i.e., $x\perp_B y$ implies $y\perp_B x$ for all $x,y\in \X$;
\smallskip
\item[(iii)] Birkhoff-James orthogonal is left-additive, i.e., $x\perp_B y$ and $z\perp_B y$ implies $(x+z)\perp_B y$ for all $x,y,z\in \X$.
\end{enumerate}
\end{thm}

On the other hand, right-additivity of Birkhoff-James orthogonality characterizes smooth points on a Banach space.

\begin{thm}[James, \cite{James 1}]\label{thm:31:10:26::10:27}
Let $\X$ be a Banach space and $x\in \X$ be arbitrary. Then $x$ is smooth if and only if Birkhoff-James orthogonality is right-additive at $x$, i.e., $x\perp_B y$ and $x\perp_B z$ imply that $x\perp_B (y+z)$ for all $y,z\in\X$.
\end{thm} 

The idea of orthogonal complement of a closed subspace in a Hilbert space was generalized to Banach spaces through the notion of $1$-complemented subspace. A closed subspace of a
Banach space is said to be \textit{orthogonally complemented} or \textit{1-complemented} \cite{R 2} if it is the range of a norm-one projection. We now recall the notion of right-complemented subspaces in terms of the Birkhoff-James orthogonality.

\medskip

\noindent
\textbf{$\bullet$ Right-complemented subspace:} Let $\X$ be a Banach space and let $\Y$ be a closed linear subspace of $\X$. Then $\Y$ is said to be \textit{right-complemented} in $\X$, if there exists a closed subspace $\widetilde{\Y}$ of $\X$ such that $\X = \Y \oplus \widetilde{\Y}$ and $\Y \perp_B \widetilde{\Y}$, i.e. $y\perp \widetilde{y}$ for all $y\in \Y$ and $\widetilde{y}\in\widetilde{\Y}$. In such case, we write $\X= \Y \bigoplus_\SR\widetilde{\Y}$ and the subspace $\widetilde{\Y}$ is said to be a \textit{right-complement} of $\Y$. An \textit{orthogonally complemented} or $1$-\textit{complemented} subspace of a Banach space is precisely a right-complemented subspace as shown in \cite[Proposition 2.7]{JPR}.

As a direct and straight-forward consequence of Theorem \ref{thm:31:10:26::10:27}, we have the following lemma.

\begin{lem}\label{lem:201}
Let $\X$ be a smooth Banach space such that $\X = (\Y_1 \oplus \Y_2) \oplus \Y_3$, for some closed subspaces $\Y_i~(1\leq i \leq 3)$. If $\Y_i \perp_B \Y_j$ for $1\leq i < j\leq 3$, then $\X= \Y_1 \oplus_\SR (\Y_2 \oplus_\SR \Y_3)$.
\end{lem}

\smallskip

\section[S\l oci\'nski-type decomposition for Wold isometries]{S\l oci\'nski-type decomposition for Wold isometries with right-complemented ranges}\label{Sec:10}

\noindent In this Section, we produce a S\l oci\'nski-type decomposition for a doubly commuting tuple of Wold isometries with right-complemented ranges acting on a reflexive, smooth and strictly convex Banach space. We begin with the definition of unilateral shift on a Banach space which was first introduced by Faulkner and Huneycutt in \cite{Faulkner Huneycutt}.

\begin{defn}\label{defn:Banach-Unilateral}
An isometry $V: \X\to \X$ is said to be a unilateral shift if there exists a subspace $\mathcal{L}\subseteq\X$ such that
\begin{itemize}
\item[(i)] For $n > m$, $V^n \mathcal{L}\perp_B V^m \mathcal{L}$.
\smallskip
\item [(ii)] $\X=\displaystyle\bigoplus_{n=0}^{\infty}V^n\mathcal{L}$,
where the direct sum $\displaystyle\bigoplus_{n=0}^{\infty}V^n\mathcal{L}$ stands for the following set:
\[
\bigoplus_{n=0}^{\infty}V^n\mathcal{L}= \{x\in \X: x \text{ can be uniquely represented as } x = \sum_{n=0}^\infty V^n l_n, ~~ l_n\in \mathcal{L} \text{ for all } n\geq 0 \}.
\]
\end{itemize}
The subspace $\mathcal{L}$ is called a \textit{generating subspace} for $V$ and the $dim(\mathcal{L})$ is called the \textit{multiplicity} of $V$. It is easy to see that the multiplicity of a unilateral shift is well-defined and uniquely determined by $\dim(\X/ V(\X))$. 
\end{defn}

Recall that for an isometry $V$ on a Banach space $\X$, an operator $\Delta \in \mathcal B(\X)$ with $\|\Delta\|=1$ is called a \textit{left-inverse} of $V$ if $\Delta V =I_{\X}$. The left-inverse obtained in \cite[Theorem 5.3]{JPR} has some additional useful properties which we record below. These properties will be used in sequel. We skip the proof as it mostly follows from the proof of Theorem 5.3 in \cite{JPR} and the rest is obvious.

\begin{prop}\label{thm:202}
Let $V$ be a Wold-isometry on a Banach space $\X$ such that $\X = \X_u \bigoplus \X_s$, $V|_{\X_u}$ is uniatry and $V|_{\X_s}$ is a pure isometry. Suppose $\mathcal{L}$ is a closed subspace such that $\X = V(\X) \oplus_\SR \mathcal{L}$ and $\X_s = \displaystyle \bigoplus_{n=0}^\infty V^n \mathcal{L}$. Then the operator $\Delta : \X \to \X$ defined by 
\begin{equation}\label{eq:204}
   \Delta ( x_u + \sum_{n=0}^\infty V^nl_n ) = (V|_{\X_u})^{-1}(x_u) + \sum_{n=0}^\infty V^nl_{n+1}, ~x_u\in \X_u, \quad l_n\in \mathcal{L}~ \forall n\geq 0,
\end{equation}
satisfies the following;
\begin{enumerate}
\item[(i)] $\Delta $ is a bounded linear map and $\|\Delta \| =1.$
\item[(ii)] Both $\X_u$ and $\X_s$ are invariant under $T.$ 
\item[(iii)] $\Delta V= I_{\X}$ and $V\Delta $ is a norm-one projection on $V(\X)$.
\item[(iv)] For every $x_s = \sum_{n=0}^\infty V^nl_n \in \X_s$, we have 
\[
  \| \Delta ^k (x_s) \| = \left\|\sum_{n=0}^\infty V^nl_{n+k} \right\| = \left\|\sum_{n=k}^\infty V^nl_n \right\| \longrightarrow 0 \quad \text{ as } k \to \infty.
\]
\end{enumerate}
\end{prop}

The operator $\Delta$ defined in (\ref{eq:204}) changes with the change of right-complement $\mathcal{L}$ of $V(\X)$. Consequently, left-inverse of a Wold isometry may not be unique. However, if the underlying Banach space is reflexive, smooth and strictly convex, then left-inverse of an isometry when exists, becomes unique as right-complement of a linear subspace of $\X$ (if exists) is unique then. Evidently, the properties $(i)-(iv)$ in Theorem \ref{thm:202} does not change with the change of $\mathcal{L}.$ In the next lemma, we show that the uniqueness of the Wold decomposition, the left-inverse and the right-complement of a Wold isometry are directly connected with each other.

\begin{lem}\label{lem:202}
Let $V$ be an isometry on a reflexive Banach space $\X$ whose range is right-complemented in $\X$. Then the following are equivalent:
\begin{enumerate}
\item[(i)] $V(\X)$ has unique right-complement in $\X$ ;
\item[(ii)] $V$ has unique left-inverse ;
\item[(iii)] $V$ has unique Wold decomposition.
\end{enumerate}
Moreover, if $\X$ smooth and strictly convex, then the assertions (i), (ii), (iii) hold .\end{lem}

\begin{proof}
The equivalence (i) and (ii) follows from (iii)-Proposition \ref{thm:202}. The implication (ii) $\Rightarrow$ (iii) follows from the definition of left-inverse of a Wold isometry and (iii)$\Rightarrow$ (i) follows from \cite[Corollary 1]{CFS}.
\end{proof}

Recall from the `Introduction' that such a tuple of isometries $(V_1, \dots ,V_n)$ is called doubly commuting if $\Delta_iV_j=V_j\Delta_i$ for $1\leq i, j \leq n$ and $i \neq j$, where $\Delta_i$ is the unique left-inverse of $V_i$ for each $i$. The assumptions of reflexivity, smoothness and strict convexity on a Banach space provide freedom of having the existence and uniqueness of left-inverse of a Wold isometry with right-complemented range acting on it, which is shown in Lemma \ref{lem:202}. Hence, we shall assume this fact (of existence and uniqueness of left-inverse) throughout this Section and the next two Sections and apply it in the proofs of the results without mentioning it. Also, the fact that a Wold isometry (see Definition \ref{Definition:Wold-isometry-003}) $V$ admits a Wold decomposition $V= V|_{\X_u} \oplus V|_{\X_s}$, where $\X_u, \X_s$ are as in $(\ref{Wold decomposition Banach spaces})$, will be used without a mention. 

\smallskip

We need a few more preparatory results to reach our desired result, viz. S\l oci\'nski-type decomposition for Banach space isometries. 

\begin{lem}\label{lem:301}
Let $(V_1, V_2)$ be a commuting pair of Wold isometries with right-complemented ranges on a Banach space. Suppose $\Delta_i$ is a left-inverse of $V_i$ for $i=1,2$. If $\Delta_iV_j = V_j \Delta_i$ for distinct $i,j\in\{1,2\}$, then $\Delta_1 \Delta_2 = \Delta_2 \Delta_1.$
\end{lem}

\begin{proof}
It follows from the hypothesis that 
\[
V_2 \Delta_1 \Delta_2 V_1 = (V_2 \Delta_1) (V_1 \Delta_2) = V_2 \Delta_2.
\]

Pre-composition by $\Delta_2$ and post-composition by $\Delta_1$, give $\Delta_1 \Delta_2 = \Delta_2 \Delta_1$. 
\end{proof}

The following lemma is a notable observation by S. Banach.

\begin{lem}\label{lem:302}
Let $\X$ be a Banach space and $\{x_n\}$ be a sequence in $\X.$ Then the sequence $\{x_n\}$ is convergent if and only if the series $\sum_{n=1}^\infty (x_{n+1} - x_n)$ is convergent.
\end{lem}

\begin{lem}\label{thm:301}
Let $V$ be a Wold isometry with right-complemented range on a Banach space $\X$ and $\Delta$ be a left-inverse of $V$ associated with the Wold decomposition $V_1|_{\X_u} \oplus V_1|_{\X_s}$. Then 
\begin{enumerate}
\item[(i)] $\X_u$ is an invariant subspace of every bounded linear map commuting with $V$.
\smallskip
\item[(ii)] $\X_s$ is an invariant subspace of every bounded linear map commuting with $\Delta$.
\end{enumerate}
\end{lem}

\begin{proof}
The subspaces  $\X_u$ and $\X_s$ are given by 
\[
 \X_u = \bigcap_{n=0}^\infty V^n(\X), \qquad \X_s = \bigoplus_{n=0}^\infty V^n \mathcal{L}, \qquad \X= V(\X) \oplus_\SR \mathcal{L},
\]

(i). Let $T\in \mathcal{B}(\X)$ be such that $TV=VT$. Then for all $n\in \mathbb{N}$, we have $V^nT = TV^n$. Therefore,
\[
   T\left(\bigcap_{n=0}^\infty V^n(\X)\right) = \bigcap_{n=0}^\infty V^n T(\X) \subseteq \bigcap_{n=0}^\infty V^n(\X) = \X_u,
\]
shows that $\X_u$ is invariant under $T$.

(ii). We show that $\X_s$ is invariant under $T$ if $\Delta T=T\Delta.$ Let $x_s\in \X_s$ be arbitrary. Then we have
\begin{equation}\label{eq:301}
 Tx_s = \sum_{k=0}^n V^k (I- V\Delta) \Delta^k (Tx_s) + V^{n+1} \Delta^{n+1}Tx_s, \quad n\in \mathbb{N}.
\end{equation}
Since $(I-V\Delta)$ is a bounded projection onto $\mathcal{L}$ (see (iii)-Proposition \ref{thm:202}) and $\X= \X_u\oplus \X_s$, it follows from (\ref{eq:301}) that 
\[
 Tx_s = \displaystyle \lim_{n\to \infty} V^n \Delta^n Tx_s + \sum_{n=0}^\infty V^n(I-V\Delta)\Delta^n Tx_s,
\]
because the series $\displaystyle \sum_{n=0}^\infty V^n(I-V\Delta)\Delta^n Tx_s$ is convergent by  Lemma \ref{lem:302}. Indeed, by (iv)-Proposition \ref{thm:202} we have
\[
 \| V^n \Delta^n Tx_s \| = \| \Delta^n Tx_s \| = \| T \Delta^n x_s \| \leq \| \Delta^n x_s \| \longrightarrow 0 \text{ as } n\to \infty,
\]
where the last equality follows from $\Delta T=T\Delta.$ Consequently, we have 
\[
 Tx_s = \sum_{n=0}^\infty V^n(I-VA)A^n Tx_s \in \X_s.
\]
This completes the proof.
\end{proof}

The proof of the following corollary is evident from Lemma \ref{thm:301}.

\begin{cor}\label{cor:301}
Let $(V_1,V_2)$ be a pair of doubly commuting Wold isometries with right-complemented ranges on a reflexive, smooth and strictly convex Banach space. If $V_1= V_1|_{\X_u} \oplus V_1|_{\X_s}$ is the Wold decomposition of $V_1$, then $\X_u$ and $\X_s$ are reducing subspaces of both $V_1$ and $V_2$.
\end{cor}

Now, we are in a position to state and prove a generalization of S\l oci\'nski's theorem in Banach space setting.

\begin{thm}\label{thm:302}
Let $(V_1,V_2)$ be a pair of doubly commuting Wold isometries with right-complemented ranges on a reflexive, smooth and strictly convex Banach space $\X$. Then there exist joint reducing subspaces $\X_{uu}, \X_{su}, \X_{us}, \X_{ss}$ of $V_1$ and $V_2$ such that
\begin{enumerate}
 \item[(i)] $\X = \X_{uu} \bigoplus \X_{us} \bigoplus \X_{su} \bigoplus \X_{ss}$;
\smallskip
  \item[(ii)] both $V_1|_{\X_{uu}}$ and $V_2|_{\X_{uu}}$ are unitary;
\smallskip
  \item[(iii)] both $V_1|_{\X_{ss}}$ and $V_2|_{\X_{ss}}$ are unilateral shifts;
\smallskip
  \item[(iv)] $V_1|_{\X_{us}}$ is unitary and $V_2|_{\X_{us}}$ is unilateral shift;
\smallskip
  \item[(v)] $V_1|_{\X_{su}}$ is unilateral shift and $V_2|_{\X_{su}}$ is unitary.
\end{enumerate}
\end{thm}

\begin{proof}
Let $\Delta_1$ and $\Delta_2$ be the left-inverse of $V_1$ and $V_2$ respectively. Let $V_1= V_1|_{\X_u} \oplus V_1|_{\X_s}$ be the Wold decomposition of $V_1$. Then $\X_u$ and $\X_s$ are reducing subspaces of both $V_1, V_2$ by Corollary \ref{cor:301}. It follows from (iii)-Proposition \ref{thm:202} that $V_2 \Delta_2|_{\X_u} : \X_u \to \X_u$ and $V_2 \Delta_2|_{\X_s} : \X_s \to \X_s$ are a norm-one projections on $V_2(\X_u)$ and $V_2(\X_s)$ respectively. Consequently, by \cite[Corollary 1]{CFS}, $V_2|_{\X_u}$ and $V_2|_{\X_s}$ are Wold isometries with right-complemented ranges. So, we have $\X_u = \X_{uu} \bigoplus \X_{us}$ and $\X_s = \X_{su} \bigoplus \X_{ss}$ such that
\begin{enumerate}
\item[(a)] $V_2 |_{\X_{uu}} : \X_{uu} \to \X_{uu}$ and $V_2 |_{\X_{su}} : \X_{su} \to \X_{su}$ are unitary (surjective isometry).
\smallskip
\item[(b)] $V_2 |_{\X_{us}} : \X_{us} \to \X_{us}$ and $V_2 |_{\X_{ss}} : \X_{ss} \to \X_{ss}$ are unilateral shifts.
\end{enumerate}
By Lemma \ref{lem:202}, the left-inverse of $V_2|_{\X_u}$ and $V_2|_{\X_s}$ is given by $\Delta_2|_{\X_u}$ and $\Delta_2|_{\X_s}$ respectively. Thus, all the subspace $\X_{uu}, \X_{us}, \X_{su}$ and $\X_{ss}$ are reducing for both of $V_1$ and $V_2$ by Corollary \ref{cor:301}. Since $\X_{uu}$, $\X_{us}\subseteq \X_u,$ we have $V_1|_{\X_{uu}}$ and $V_1 |_{\X_{us}}$ are unitary. On the other hand, since $\X_{su}$, $\X_{ss}\subseteq \X_s$, $V_1|_{\X_{ss}}$ and $V_1 |_{\X_{ss}}$ are unilateral shifts. Consequently, $\X = \X_{uu} \bigoplus \X_{us} \bigoplus \X_{su} \bigoplus \X_{ss}$, where $\X_{uu}, \X_{us}, \X_{su}$ and $\X_{ss}$ are the subspaces satisfying the stated properties.
\end{proof}

The above theorem can be generalized from a pair to a tuple of such Banach space isometries following the same algorithm. Only a little more elaborative computations can occur. For this reason, we present this result below without a proof. We follow the standard convention that for any nonempty sets $A$ and $B$, $B^A$ denotes the set of all functions from $B$ to $A$.

\begin{thm}\label{thm:303}
Let $(V_1, \dotsc, V_m)$ be a tuple of doubly commuting Wold isometries with right-complemented ranges on a reflexive, smooth and strictly convex Banach space $\X$. Let $A=\{u,s\}$ and $B=\{1,2,\dotsc, m\}$. Then there are $2^m$ many subspaces $\left\{\X_f: f\in B^A\right\}$ such that the following holds:
\begin{enumerate}
\item[(i)] $\X_f$ is a joint reducing subspace of $\{V_i\}_{_1}^m$ for each $f\in B^A$;

\item[(ii)] for $1\leq i \leq m$, and $f\in B^A$, $V_i|_{\X_f}$ is unitary if $f(i) =u$ and unilateral shift if $f(i) =s$; 

\item[(iii)] for the constant function $f(j)=u$, where $(f\in B^A)$, $\X_f$ is the maximal joint reducing subspace of $\{V_i\}_{_1} ^m$ such that $V_{i}|_{\X_f}$ is a unitary for $1\leq i \leq m$;

\item[(iv)] for the constant function $g(j)=s$, where $(g\in B^A)$, $\X_g$ is the maximal joint reducing subspace of $\{V_i\}_{_1} ^m$ such that $V_{i}|_{\X_g}$ is a unilateral shift for $1\leq i \leq m$.
\end{enumerate}
\end{thm}

\smallskip

\section[A Berger-Coburn-Lebow model for commuting Wold isometries]{A Berger-Coburn-Lebow model for commuting Wold isometries with right-complemented ranges}\label{Sec:11}

\noindent
Throughout this Section, we consider isometries that are Wold isometries with right-complemented ranges acting on reflexive, smooth and strictly convex Banach spaces. As we have mentioned in the beginning of Section \ref{Sec:10}, these three assumptions on a Banach space guarantees (by Lemma \ref{lem:202}) the existence and uniqueness of left-inverse of a Wold isometry with right-complemented range acting on it. We shall continue applying this fact (existence and uniqueness of left-inverse) in the proofs of the results without mentioning it. The main aim of this Section is to present a structure-cum-model theorem for a commuting tuple of such Wold isometries. This generalizes the renowned Berger-Coburn-Lebow (BCL) model (e.g., see \cite{BCL}) for a tuple of commuting Hilbert space isometries, which is stated below.

\begin{thm}[Berger-Coburn-Lebow, \cite{BCL}]\label{thm:103} Let $(V_1, \dotsc, V_m)$ be a tuple of commuting isometries on a Hilbert space $\HS$ such that $V= \prod_{i=1}^m V_i$ is a pure isometry. Then there exist tuples $(P_1, \dotsc, P_m)$ and $(U_1, \dotsc, U_m)$ consisting of orthogonal projections and unitaries in $\mathcal{B}(\mathcal{D}_{V^*})$, respectively, such that
\[
   (V_1, \dotsc, V_m, V) \equiv (M_{P_1^{\perp} U_1 + z P_1 U_1}, \dotsc, M_{P_m^{\perp} U_m + z P_m U_m}, M_z),,
\]
on $H^2(\mathcal{D}_{V^*})$ where $\mathcal{D}_{V^*} = (I - VV^*)(\mathcal{H})$ and $ P_i^\perp = (I - P_i)$ for all $1\leq m$.
\end{thm}

An independent proof to the above theorem on BCL model is also given in \cite{P & S}. A complete classification of a
commuting $m$-tuple $(V_1, \dotsc, V_m)$ with their product $V = \prod_{i=1}^m V_i$ being a pure isometry of multiplicity $m$, is established in \cite{BDF}. Also, in \cite{BDF II} a canonical model is obtained for a commuting pair $V=(V_1, V_2)$ of Hilbert space isometries by associating with $V$ a unilateral shift $V_0$ and a contraction $T$ that acts as a commutant of $V_0$. Commutants of certain kind of the unilateral shift $M_z$ play the key role in the BCL model for commuting Hilbert space isometries. For a Hilbert space $\HS$, a commutant of the unilateral shift $M_z$ on $H^2(\HS)$ is isometrically isomorphic (i.e., unitarily equivalent) with $M_{\varphi}$ for some $\varphi \in H^\infty(\mathcal{B}(\HS))$. For the convenience of a reader, let us recall that the Hardy space $H^2(\mathcal H)$ consists of analytic functions $f:\D \rightarrow \mathcal H$ that have square summable coefficients in their power series, i.e.,
\begin{equation*}
H^2(\mathcal H)=\{ f: \D \rightarrow \mathcal H \, : f \text{ is analytic, } f(z)=\sum_{n=0}^{\infty} a_nz^n \,, \, a_n \in \mathcal H \, \& \, \sum_{n=0}^{\infty}\, \|a_n\|^2 < \infty \}
\end{equation*}
and the multiplication operator $M_{\varphi}: H^2(\mathcal H) \rightarrow H^2(\mathcal H)$ for $\varphi \in H^{\infty}\left( \mathcal B(\mathcal H) \right)$ is defined by
\[
M_{\varphi}(f)(z):=\varphi(z)\left( f(z) \right).
\]
The commutant of a Banach space shift of multiplicity one was found by Crownover in \cite{RMC}, though the same for unilateral shift of arbitrary multiplicity is not known till date (as per our knowledge). Thus, as a first step towards a possible BCL model for Banach space isometries is to determine the commutants of a Banach space unilateral shift with right-complemented range and of arbitrary multiplicity. Note that the range of a unilateral shift may not be right-complemented in general, e.g., see \cite[Example 5.2]{JPR}. However, the unilateral shift part of a Wold isometry with right-complemented range also has right-complemented range. For determining the commutants of a Banach space unilateral shift with right-complemented range and of arbitrary multiplicity, we need to recall a few elementary concepts.

\smallskip

\noindent \textbf{$\bullet$ Sequence space associated with a unilateral shift on a Banach space:} Let $V$ be a unilateral shift on a Banach space $\X$ with generating subspace $\mathcal{L}$. Consider the sequence space
\begin{equation}\label{eq:201}
  \mathcal{S} =\left\{\{l_n\}_{n\in \N_0}: l_n\in \mathcal{L} \text{ for all } n\in \mathbb{N}\cup\{0\}, ~ \sum_{n=0}^\infty V^n l_n \text{ converges in } \X\right\}, 
\end{equation}
and is equipped with the norm
\[
  \left\|\{l_n\}_{n\in \N_0}\right\|_\mathcal{L}= \left\| \sum_{n=0}^\infty V^n l_n \right\|, \quad \{l_n\}_{n\in \N_0}\in \mathcal{S}.
\]
Consider $\Lambda: \X \to \mathcal{S}$ defined by $\Lambda \left(\sum_{n=0}^\infty V^n l_n\right) = \{l_n\}_{n\in \N_0}$. Then $\Lambda$ is a unitary and $\X$ is isometrically isomorphic with $\mathcal{S}$ via $\Lambda$. The isometry $V$ can be realised as the forward shift operator $\widehat{V}$ on $\mathcal{S}$ which is given by $\widehat{V} = \Lambda V \Lambda^{-1}$, i.e.,
\begin{equation}\label{eq:202}
 \Lambda V \Lambda^{-1}(\{l_n\}_{n\in \N_0}) = \widehat{V}(\{l_n\}_{n\in \N_0}) = \{l'_n\}_{n\in \N_0}, \quad l'_n = \begin{cases}
  \mathbf{0} \quad \text{ if } n=0, \\
  l_n \quad \text{ if } n\geq 1.
  \end{cases}
\end{equation}
Then $\widehat{V}$ is also a unilateral shift with generating subspace $\widehat{L}=\{(l,\mathbf{0}, \mathbf{0}, \dotsc):l\in \mathcal{L}\}$.

\medskip

The notations $\mathbf{T_L}\left[S_{i-j}\right]$, $\mathbf{T}_{A+zB}$ and $\mathbf{T}_{A+\bar{z}B}$ are used to represent the following Toeplitz block matrices on $\mathcal{S}$ :
\begin{equation}\label{eq:203}
\mathbf{T_L}\left[S_{i-j}\right] = \begin{bmatrix}
          S_0 & 0 & 0 & \cdots \\
          S_1 & S_0 & 0 & \cdots \\
          S_2 & S_1 & S_0 & \cdots \\
          \vdots & \vdots & \vdots & \ddots
         \end{bmatrix}, \ \mathbf{T}_{A+zB} = \begin{bmatrix}
          A & 0 & 0 & \cdots \\
          B & A & 0 & \cdots \\
          0 & B & A & \cdots \\
          \vdots & \vdots & \vdots & \ddots
         \end{bmatrix}, \ \mathbf{T}_{A+ \bar{z}B} = \begin{bmatrix}
          A & B & 0 & \cdots \\
          0 & A & B & \cdots \\
          0 & 0 & A & \cdots \\
          \vdots & \vdots & \vdots & \ddots
         \end{bmatrix},
\end{equation}
where all of $A, B$ and $S_i$ ($0\leq i <\infty$) are bounded operators on $\mathcal{L}$. Now, we state the theorem that determines the commutants of a Banach space unilateral shift with right-complemented range and of arbitrary multiplicity.

\begin{thm}\label{thm:401}
Let $V$ be a unilateral shift on a Banach space $\X$ with generating subspace $\mathcal{L}$ and $\X = V(\X) \oplus_\SR \mathcal{L}$. Consider the forward shift operator $\widehat{V}\in \mathcal{B}(\mathcal{S})$ with generating subspace $\widehat{\mathcal{L}}$, where $(\mathcal{S}, \|\cdot\|_\mathcal{L})$ and $\widehat{V}$ are defined as in $(\ref{eq:201})$ and $(\ref{eq:202})$. Then $ T\in \mathcal{B}(\mathcal{S})$ commutes with $\widehat{V}$ if and only if there exists a uniformly bounded sequence $\{T_n\}_{n\in \N_0}$ such that $T_n\in \mathcal{B}(\mathcal{L})$ for all $n\geq 0$ and the block matrix representation of $T$ with respect to the decomposition $\mathcal{S}= \bigoplus_{n=0}^\infty \widehat{V}^n \widehat{\mathcal{L}}$ is a Toeplitz block matrix given by $\mathbf{T_L}\left[T_{i-j}\right]$.
\end{thm}

\begin{proof}
First we show that, an operator $T\in \mathcal{B}(\mathcal{S})$ with block matrix representation as in the statement commutes with $\widehat{V}.$ Let $\{l_n\}_{n\in \N_0}\in \mathcal{S}$ be arbitrary. Then 
\begin{align*}
    & \widehat{V} ~ T (\{l_n\}_{n\in \N_0}) = (\mathbf{0}, a_0, a_1, a_2, \cdots), \qquad a_n = \sum_{k=0}^n T_{n-k}l_k, \quad n\geq 0 \quad \text{ and } \\
    & T ~ \widehat{V}(\{l_n\}_{n\in \N_0}) = T (\mathbf{0}, l_0, l_1, l_1, \cdots) = (\mathbf{0}, a_0, a_1, a_2, \cdots) = \widehat{V} ~ T (\{l_n\}_{n\in \N_0}),
\end{align*}
shows that $T\widehat{V} = \widehat{V} T.$ For the other direction, suppose $T\in \mathcal{B}(\mathcal{S})$ such that $T\widehat{V} = \widehat{V} T.$ Since $\X = V(\X) \oplus_\SR \mathcal{L}$, we have $\mathcal{S} = \widehat{V}(\mathcal{S}) \oplus_\SR \widehat{\mathcal{L}}$. Let $\widehat{\Delta}$ be a left-inverse of $\widehat{V}.$ For every $n\geq 0,$ define $T_n: \mathcal{L} \to \mathcal{L}$ by 
\[
    T_n (l) = \Lambda^{-1} (I-\widehat{V} \widehat{\Delta})\widehat{\Delta}^n T \Lambda(l),\quad l\in \mathcal{L},
\]
where $\Lambda : \X \to \mathcal{S}$ is the unitary as in (\ref{eq:202}). Then for every $n\geq 0,$ $T_n\in \mathcal{B}(\mathcal{L})$ with $\|T_n\| \leq \|(I-\widehat{V}\widehat{\Delta})\widehat{\Delta}^n \; T\| \leq 2 \|T\|.$ It is easy easy to see that 
\[
   T(l, \mathbf{0}, \cdots) = (T_0(l), T_1(l), \cdots), \qquad l\in \mathcal{L}.
\]
The subspace $C_{00}(\mathcal{L})$ consisting of all sequences with zero entries except for finitely many is dense in $\mathcal{S}$. Indeed, for any $\{l_n\}_{n\in \N_0}\in \mathcal{S},$ we have
\[
    \left\| (l_0, \cdots, l_m, \mathbf{0}, \cdots) - \{l_n\}_{n\in \N_0} \right\| = \left\| \sum_{k=m+1}^\infty V^k l_k \right\| \longrightarrow 0 \quad \text{ as } m\to \infty.
\]
Consider the operator $\widehat{F}: C_{00}(\mathcal{L}) \to \mathcal{S}$ defined by 
\[
    \widehat{F}(l_0, \cdots, l_n, \mathbf{0}, \cdots) = (x_0, x_1, \cdots), \quad x_k = \begin{cases}
                     \sum_{i=0}^k T_{k-i}(l_i) & \text{ if } 0\leq k \leq n , \\
                     \sum_{i=0}^n T_{n+r-i}(l_i) & \text{ if } k = n+ r, ~ r\geq 1. 
                     \end{cases}
\]
Then for every $n\geq 1$ and $(l_0, \cdots, l_n, \mathbf{0}, \cdots) \in C_{00}(\mathcal{L})$, we have
\begin{align*}
   T(l_0, \cdots, l_n, \mathbf{0}, \cdots) = T\left(\sum_{i=0}^n \widehat{V}^i (l_i, \mathbf{0}, \cdots) \right) & = \sum_{i=0}^n \widehat{V}^i T(l_i, \mathbf{0}, \cdots),\left[\text{ as }\widehat{V}T = T\widehat{V} ~\right] \\
   & = \sum_{i=0}^n \widehat{V}^i (T_0l_i ,T_1 l_i, \cdots) \\
   & = \widehat{F}(l_0, \cdots, l_n, \mathbf{0}, \cdots).
\end{align*}
This shows that $\widehat{F}$ is a well-defined linear map and $\|\widehat{F}\| \leq \|T\|.$ Therefore, $\widehat{F}$ uniquely extends to a bounded linear map $F: \mathcal{S} \to \mathcal{S}$. More precisely, for every $\eta=\{l_n\}_{n\in \N_0}\in \mathcal{S}$, $F(\eta) = \displaystyle \lim_{n\to \infty} \widehat{F}(\eta_n)$, where $\eta_n = (l_0, \dotsc, l_n, \mathbf{0}, \dotsc)$ for all $n\geq 0$. In addition, we have
\begin{align*}
  F(\eta) = \lim_{n\to \infty} \widehat{F}(\eta_n) = \lim_{n\to \infty} T(\eta_n )  = T(\eta) & = T \sum_{n=0}^\infty \widehat{V}^n (l_n,\mathbf{0},\cdots) \\
  & = \sum_{n=0}^\infty \widehat{V}^n T(l_n,\mathbf{0}, \cdots)~~\left[\text{ as } \widehat{V} T = T \widehat{V} ~ \right] \\
  & = \sum_{n=0}^\infty \widehat{V}^n \sum_{k=0}^\infty \widehat{V}^k (T_k l_n, \mathbf{0}, \cdots) \\
  & = \sum_{n=0}^\infty \sum_{k=0}^\infty \widehat{V}^{n+k} (T_k l_n, \mathbf{0}, \cdots) \\
  & = \sum_{m=0}^\infty \widehat{V}^m (a_m, \mathbf{0}, \cdots), ~~ a_m = \sum_{i+j =m}T_il_j,~ m\geq 0.
\end{align*}
Consequently, the block matrix representation of the operator $F=T$ is a Toeplitz block matrix and is given by $\mathbf{T_L}\left[T_{i-j}\right]$. This completes the proof.
\end{proof}

\begin{note}
The operator $T\in \mathcal{B}(\mathcal{S})$ commuting with $\widehat{V}$ is uniquely determined by its image of the subspace $\{T(l,\mathbf{0}, \cdots): l\in \mathcal{L}\}$.
\end{note}

\begin{cor}\label{cor:401}
Suppose $V$ is a unilateral shift operator on a Banach space $\X$ with generating subspace $\mathcal{L}.$ Then $\mathbf{T}\in \mathcal{B}(\X)$ commutes with $V$ if and only if there exists a uniformly bounded sequence $\{T_n\}_{n\in \N_0}$ of operators on $\mathcal{L}$ such that the Toeplitz block matrix representation of $\mathbf{T}$ is given by $\Lambda^{-1}\mathbf{T_L}\left[T_{i-j}\right]\Lambda$, where $\Lambda: \X \to \mathcal{S}$ is the unitary as in $(\ref{eq:202})$. Moreover,
\[
   \mathbf{T}(l) = \sum_{n=0}^\infty V^n T_n l, \qquad l\in \mathcal{L}. 
\]
\end{cor}

\subsection{Berger-Coburn-Lebow (BCL) representation}\label{Subsec:41} 

In this Subsection, we present a BCL model for a tuple of commuting Wold isometries with right-complemented ranges acting on a reflexive, smooth and strictly convex Banach space. We need to have a few preparatory results before that.

\begin{lem}\label{lem:411}
Let $(V_1,V_2)$ be a pair of commuting Wold isometries with right-complemented ranges acting on a reflexive, smooth and strictly convex Banach space $\X$. Suppose $\mathcal{L}_i$ is the right-complement of $V_i(\X)$ in $\X$ and $\Delta_i$ is the left-inverse of $V_i$ for $i=1,2$. Then
\begin{enumerate}
\item[(i)] $V = V_1V_2 = V_2V_1$ is a Wold isometry with right-complemented range and the right-complement $\W$ of $V(\X)$ is given by
\[
    V_1\mathcal{L}_2 \oplus_\SR \mathcal{L}_1 = \W = V_2\mathcal{L}_1 \oplus_\SR \mathcal{L}_2.
\]

\smallskip

\item[(ii)] The left-inverse $\Delta$ of $V$ is given by $\Delta = \Delta_1 \Delta_2 = \Delta_2\Delta_1$.

\smallskip

\item[(iii)] The subspace $\W$ is invariant under $\Delta_i$ $(i=1,2)$ and $\ker(\Delta) = \W$.

\item[(iv)] There are norm-one projections $Q_1$ and $Q_2$ on $\W$ such that every $w\in \W$ can be uniquely expressed as 
\[
Q_1(w) + (I_\W - Q_1)(w) = w = Q_2(w) + (I_\W - Q_2)(w).
\]
\end{enumerate}
\end{lem}

\begin{proof}
(i) Since $V_1,V_2$ are linear isometries, they preserve Birkhoff-James orthogonality. Thus, it follows from the hypothesis that $V_1(\X) = V_1 V_2(\X) \bigoplus_\mathsf{R} V_1\mathcal{L}_2$ and $V_2(\X) = V_2 V_1(\X) \bigoplus_\mathsf{R} V_2\mathcal{L}_1.$ Thus, we have
\[
  \left( V_1 V_2(\X) \oplus_\SR V_1\mathcal{L}_2 \right) \oplus_\SR \mathcal{L}_1 = \X = \left( V_2 V_1(\X) \oplus_\SR V_2\mathcal{L}_1 \right) \oplus_\SR \mathcal{L}_2.
\]
Since $\X$ is smooth and $V_1 V_2 = V_2 V_1$, by Lemma \ref{lem:201} we have 
\[
 V_1V_2(\X) \oplus_\SR \left(V_1\mathcal{L}_2 \oplus_\SR \mathcal{L}_1\right) = \X = V_1 V_2 (\X) \oplus_\SR\left(V_2\mathcal{L}_1 \oplus_\SR \mathcal{L}_2\right).
\]
Therefore, $V =V_1 V_2 = V_2V_1$ is a Wold isometry with right-complemented range by Corollary 1 in \cite{CFS}. By virtue of Lemma \ref{lem:202}, the right-complement $\W$ of $V(\X)$ is given as stated.

\smallskip

(ii) Since $V_1V_2 = V_2 V_1$, we have $\Delta_1\Delta_2 V = I_{\X} = \Delta_2\Delta_1 V$. Therefore, by Lemma \ref{lem:202}, we have $\Delta = \Delta_1\Delta_2 = \Delta_2 \Delta_1$.

\smallskip

(iii) Note that $\Delta_i(\W) = \Delta_i \left(V_i\mathcal{L}_j \oplus\mathcal{L}_i\right) = \mathcal{L}_j \subseteq V_j \mathcal{L}_i \oplus \mathcal{L}_j = \W$ for $i,j\in \{1,2\}$ and $i\neq j$. Hence $\W$ is invariant under $\Delta_i$ for $i=1,2$. It follows from the definition of left-inverse of a Wold isometry with right-complemented range (see (\ref{eq:204})) that $\ker(\Delta) = \W$. 

\smallskip

(iv) Let $w\in \W$ be arbitrary. Then there exist unique vectors $\{l_1,r_1\}$ in $\mathcal{L}_1$ and $\{l_2, r_2\}$ in $\mathcal{L}_2$ such that $V_1l_2 + l_1 = w = V_2r_1 + r_2$. Define $Q_i: \W \to \W$ by 
\[
 Q_1(V_1l_2 + l_1)= V_1\Delta_1(w) = V_1l_2 \quad \text{ and } \quad Q_2(V_2 r_1 + r_2) = V_2\Delta_2(w) = V_2 r_1.   
\]
Then by virtue of $V_i(\X) \perp_B \mathcal{L}_i$ ($i=1,2$) we have
\[
 \|Q_1(V_1l_2+l_1)\| = \|V_1 l_2 \| \leq \|V_1l_2 +l_1\|, \quad \|Q_2( V_2r_1+r_2)\| = \|V_2r_1\| \leq \|V_2r_1 +r_2\|,
\]
and $Q_i^2 = Q_i$ for $i=1,2$. Thus, $Q_i$ is a norm-one projection for $i=1,2$. Consequently, 
\[
 Q_1(w) + (I_\W - Q_1)(w) = w = Q_2(w) + (I_\W - Q_2)(w).
\]
This completes the proof.
\end{proof}

\begin{lem}\label{lem:412}
Let $(V_1, \dotsc, V_m)$ be a tuple of commuting isometries acting on a Banach space $\X$. If $V = \displaystyle \Pi_{i=1}^m V_i$ is a unitary $($surjective isometry$)$, then $V_i$ is unitary for all $1\leq i \leq m.$
\end{lem}
\begin{proof}
It is enough to prove for a pair $(V_1, V_2)$ of commuting isometries. We show that both $V_1$ and $V_2$ are surjective. Let $y\in \X$ be arbitrary. Then there is $x\in \X$ such that $V_1 V_2x = V_2 V_1x = y$. Consequently, $V_1x$ and $V_2x$ are the inverse images of $y$ under $V_2$ and $V_1$ respectively.
\end{proof}

\begin{prop}\label{prop:411}
Let $(V_1, \dotsc, V_m)$ be a tuple of commuting Wold isometries with right-complemented ranges on a reflexive, smooth and strictly convex Banach space $\X$. Suppose $\Delta_i$ is the left-inverse of $V_i$ and $\mathcal{L}_i$ is the right-complement of $V_i(\X)$ for $1\leq i \leq m$. Let $S_m$ be the set of all permutations on $\{1,2, \dotsc, m\}$. Then the following hold.
\begin{enumerate}
\item[(i)] $V=\prod_{i=1}^m V_i$ is a Wold isometry with right-complemented range.

\smallskip

\item[(ii)] Let $\W$ be the right-complement of $V(\X)$. Then for every $\rho \in S_m$, $\W$ can be expressed as 
\begin{equation}\label{eq:401}
\left\{\begin{split}
& \W = \bigoplus_{i=1}^{m-1} \left( \prod_{j=1}^i V_{\rho(j)} \mathcal{L}_{\rho(i+1)} \right) \oplus \mathcal{L}_{\rho(1)}, \quad \text{ and} \\
&\left(\prod_{i=1}^k V_{\rho(i)}\mathcal{L}_{\rho(k+1)}\right) \perp_{B} \left(\bigoplus_{i=1}^{k-1} \prod_{j=1}^i V_{\rho(j)}\mathcal{L}_{\rho(i+1)} \oplus \mathcal{L}_{\rho(1)} \right), 1\leq k \leq m-1,
\end{split}\right.
\end{equation}
where direct sum over empty set is excluded.
\smallskip

\item[(iii)] The left-inverse $\boldsymbol{\Delta}$ of $V$ is given by $\boldsymbol{\Delta} = \displaystyle \prod_{i=1}^m \Delta_i$ and is independent of ordering. \\

\smallskip

\item[(iv)] $\W$ is invariant under $\Delta_i$ for $1\leq i \leq n$ and $\ker(\boldsymbol{\Delta}) = \W$.

\smallskip

\item[(v)] Let $\X = \X_u \oplus \X_s$ be the Wold decomposition of $V$. Then $\X_s$ is invariant under $V_i$ and $V_i|_{\X_s}$ is a Wold isometry with right-complemented range for $1\leq i \leq m$.

\end{enumerate}
\end{prop}
\begin{proof}
(i) Follows from repeated application of (i)-Lemma \ref{lem:411}.

\smallskip

(ii) We apply mathematical induction on the number of isometries. The result is true for $m=2$ by Lemma \ref{lem:411}. Suppose the result is true for any tuple of $m-1$ ($m\geq 3$) isometries and $\rho \in S_m$. Then by induction hypothesis, the right-complement $\widetilde{\W}$ of $\prod_{i=2}^m V_{\rho(i)}$ can be represented as
\begin{align*}
\begin{gathered}
 \widetilde{\W} = \bigoplus_{i=2}^{m-1} \left(\prod_{j=2}^i V_{\rho(j)} \mathcal{L}_{\rho (i+1)} \right) \oplus \mathcal{L}_{\rho(2)}, \text{ and } \\
 \left(\prod_{i=2}^k V_{\rho(i)}\mathcal{L}_{\rho(k+1)}\right) \perp_{B} \left(\bigoplus_{i=2}^{k-1} \prod_{j=2}^i V_{\rho(j)}\mathcal{L}_{\rho(i+1)} \oplus \mathcal{L}_{\rho(2)} \right), \ \ 2\leq k \leq m-1,
 \end{gathered}
\end{align*}
where direct sum over empty set is excluded. Since $V_{\rho(1)}$ is an isometry, it preserves Birkhoff-James orthogonality. Therefore, we have 
\begin{equation}\label{eq:402}
\left\{\begin{split}
& V_{\rho(1)}(\widetilde{\W}) = \bigoplus_{i=1}^{m-1} \prod_{j=1}^i V_{\rho(j)} \mathcal{L}_{\rho(i+1)}, \text{ and } \\
&\left(\prod_{i=1}^k V_{\rho(i)} \mathcal{L}_{\rho(k+1)}\right) \perp_B \left( \bigoplus_{i=1}^{k-1} \prod_{j=1}^i V_{\rho(j)}\mathcal{L}_{i+1}\right),  \ \ 2 \leq k \leq m-1 .
\end{split}\right.
\end{equation}
Now, $\X= \displaystyle\prod_{i=2}^m V_{\rho(i)} (\X) \oplus_\SR \widetilde{\W}$ and $\X = V_{\rho(1)} (\X) \oplus_\SR \mathcal{L}_{\rho(1)}$ implies that
\[
\X = \left( \prod_{i=1}^m V_{\rho(i)}(\X) \oplus_\SR V_{\rho(1)}(\widetilde{\W})\right) \oplus_\SR \mathcal{L}_{\rho(1)}.
\]
Therefore, by Lemma \ref{lem:201} we have $\X = V(\X) \oplus_\SR \left(V_{\rho(1)}(\widetilde{\W})\oplus_\SR V_{\rho(1)} \right)$. Also, by (\ref{eq:402}) we have
\begin{align*}
\begin{gathered}
   V_{\rho(1)}(\widetilde{\W}) \oplus_\SR \mathcal{L}_{\rho(1)} = \bigoplus_{i=1}^{m-1} \left( \prod_{j=1}^i V_{\rho(j)} \mathcal{L}_{\rho(i+1)} \right) \oplus \mathcal{L}_{\rho(1)}, \quad \text{ and} \\
   \left(\prod_{i=1}^k V_{\rho(i)}\mathcal{L}_{\rho(k+1)}\right) \perp_{B} \left(\bigoplus_{i=1}^{k-1} \prod_{j=1}^i V_{\rho(j)}\mathcal{L}_{\rho(i+1)} \oplus \mathcal{L}_{\rho(1)} \right), \ \ 1\leq k \leq m-1,
\end{gathered}
\end{align*}
where direct sum over empty set is excluded. By Lemma 2 in \cite{Calvert}, $V(\X)$ has unique right-complement in $\X$. Hence, the right-complement $\W$ of $V(\X)$ has the form as stated in (\ref{eq:401}) for every $\rho \in S_m$.

\smallskip

(iii) Since $(V_1, \dotsc, V_m)$ is a commuting tuple, $V= \displaystyle\prod_{i=1}^m V_{\rho(i)}$ for every permutation $\rho\in S_m$. Also, $\displaystyle \prod_{i=1}^m \Delta_{\rho(i)} V = I_\X$ holds for every $\rho\in S_m$. By Lemma \ref{lem:202}, $V$ has unique left-inverse. Consequently, the left-inverse $\boldsymbol{\Delta}$ of $V$ is given by $\displaystyle\prod_{i=1}^m \Delta_i$ and is independent of ordering. The representation of $\W$ is independent of choice of the permutation $\rho\in S_m$. For all calculations related to $V_k$ ($1\leq k \leq m$) in the proof of (iv) and (v), we use the expression of $\W$ with respect to the permutation $\sigma^k$ where $\sigma(i) = i+1$ for $1\leq i \leq m-1$ and $\sigma(m)= 1$. It is easy to see that $\sigma^m(i)=i $ for $1\leq i \leq m$. For $1\leq k < m$, $\sigma ^k(i) = i+k$ for $1\leq i \leq m-k$ and $\sigma^k (m-k+i) = i$ for $1\leq i \leq k$. 

\smallskip

(iv) Let $k\in\{1,2,\dotsc,m\}$ be arbitrary. To show $\W$ is invariant under $\Delta_k$, consider the expression of $\W$ corresponding to the permutation $\eta = \sigma^{k-1}$ if $2\leq k \leq m$ and $\eta = \sigma^m$ if $k=1$.
Then we have,
\begin{align*}
  \Delta_k (\W) = \Delta_{\eta(1)}\left(\W\right) & = \bigoplus_{i=2}^{m-1} \left(\prod_{j=2}^i V_{\eta(j)} \mathcal{L}_{\eta (i+1)} \right)\oplus \mathcal{L}_{\eta(2)} \\
  & \subseteq \W = \bigoplus_{i=1}^{m-1} \left( \prod_{j=1}^i V_{\sigma^k (j)} \mathcal{L}_{\sigma^k(i+1)} \right) \oplus \mathcal{L}_{\sigma^k (1)}.
\end{align*}
Consequently, $\W$ is invariant under $\Delta_k$ for all $1\leq k \leq m$. Also, it follows from the definition of $\boldsymbol{\Delta}$ that $\ker(\boldsymbol{\Delta}) = \W$.

\smallskip

(v) Let $k\in\{1,2,\dotsc,m\}$ be arbitrary. Then we have
\begin{align*}
   V_k (\W) = V_{\sigma^k(m)}\left(\W\right) & = \left(\prod_{i=1}^m V_{\sigma^k(i)} \mathcal{L}_{\sigma^k(m)}\right) \bigoplus \left(\bigoplus_{i=1}^{m-1} \left( \prod_{j=1}^i V_{\sigma^{k-1} (j)} \mathcal{L}_{\sigma^{k-1}(i+1)} \right) \right) \nonumber \\
   & = V\left(\mathcal{L}_k\right) \oplus \ \bigoplus_{i=1}^{m-1} \prod_{j=1}^i V_{\sigma^{k-1} (j)} \mathcal{L}_{\sigma^{k-1}(i+1)}.
\end{align*}
Consequently, we have $\X_s = V_k(\X_s) \oplus_\SR \mathcal{L}_{\sigma^{k-1}(1)} = V_k(\X_s) \oplus_\SR \mathcal{L}_{k}$, for $1\leq k \leq m$. This shows that $\X_s$ is invariant under $V_k$ for $1\leq k \leq m$. Also, Corollary 1 in \cite{CFS} implies that $V_{k}|_{\X_s}$ is a Wold isometry with right-complemented range. This completes the proof.
\end{proof}

We follow the same notations as in Proposition \ref{prop:411}. Let $w\in \W$ be arbitrary. Then for every $1\leq k\leq m$, there is an $m$ tuple $(l^{(k)}_1,\dotsc,l^{(k)}_m)$ such that $l^{(k)}_i\in \mathcal{L}_i$ for all $1\leq i \leq m$ and 
\[
   w = \sum_{i=1}^{m-1} \prod_{j=1}^i V_{\sigma^k(j)} \left( l_{\sigma^k(i+1)}^{(k)} \right) + l_{\sigma^k(1)}^{(k)}, \ 1\leq k \leq m.
\]
Since $\ker (\Delta_i) = \mathcal{L}_i$ for $1\leq i \leq m$, we have
\begin{equation}\label{eq:403}
  l^{(k)}_k = l^{(k)}_{\sigma^k(m)} = \prod_{i=1}^{m-1} \Delta_{\sigma^k(m-i)}(w) = \prod_{i=1}^{m-1} \Delta_{\sigma^k(i)} , \quad 1\leq  k \leq m,
\end{equation}
where the last equality follows from (iii)-Proposition \ref{prop:411}. For $1\leq k \leq m$, $\displaystyle\prod_{i=1}^{m-1}V_{\sigma^k(i)}$ is a Wold isometry with right-complemented range and its left-inverse is given by $\displaystyle\prod_{i=1}^{m-1} \Delta_{\sigma^k(i)}$ (by (i) \& (iii)-Proposition \ref{prop:411}). Therefore, by (iii)-Proposition \ref{thm:202}, we have 
\[
\mathbf{P_k} = \displaystyle \prod_{i=1}^{m-1} V_{\sigma^k(i)} \prod_{i=1}^{m-1} \Delta_{\sigma^k(i)}
\]
is a norm-one projection on $\W$. Consequently, we have
\begin{align}\label{eq:404}
   V_k (w) = V_{\sigma^k(m)}(w) & = \left(\prod_{i=1}^m V_{\sigma^k(i)} l^{(k)}_{\sigma^k(m)} \right) + \sum_{i=1}^{m-1} \left(\prod_{j=1}^i V_{\sigma^{k-1}(j)} l_{\sigma^{k-1}(i+1)}^{(k)} \right) \nonumber \\
   & = V \left(\prod_{i=1}^{m-1} \Delta_{\sigma^k(i)}(w)\right) + x_k, \ \ 1\leq k \leq m \ \ \left[\text{ by } (\ref{eq:403}) \right],
\end{align}
where $x_k\in \W$ for $1\leq k \leq m$ and is given by
\begin{equation}\label{eq:405}
 x_k = V_k (w) - V\prod_{i=1}^{m-1} \Delta_{\sigma^k(i)}(w) = V_k (I_{\W} -\mathbf{P_k})(w).
\end{equation}
Note that $(I_\W - \mathbf{P_k})$ is a bounded projection (may not have norm 1) and its range is given by
\begin{equation}\label{eq:406}
 (I_\W - \mathbf{P_k})(\W) = \bigoplus_{i=1}^{m-2} \prod_{j=1}^i V_{\sigma^k(j)}  \mathcal{L}_{\sigma^k(i+1)} \oplus \mathcal{L}_{\sigma^k(1)},
\end{equation}
where direct sum over empty set is excluded. Now, consider the Wold decomposition $V|_{\X_u}\oplus V|_{\X_s}$ of $V$ and the unitary $\Lambda: \X_s \to \mathcal{S}$ as in (\ref{eq:202}). By (v)-Proposition \ref{prop:411}, we have $\X_s$ is invariant under $V_k$ and $V_k|_{\X_s}$ is a Wold isometry for $1\leq k \leq m$. Also, $V_k|_{\X_s}$ commutes with the unilateral shift $V|_{\X_s}$. So, Corollary \ref{cor:401} implies that $V_k|_{\X_s}$ is a Toeplitz block matrix for every $k$. By using (\ref{eq:404}) \& (\ref{eq:405}), we have
\begin{equation}\label{eq:407}
\left\{\begin{split}
 V_k|_{\X_s}= \Lambda^{-1} \mathbf{T}_{M_k + zN_k} \Lambda~ ; ~
 & M_k = V_k (I_\W - \mathbf{P_k}), \quad \text{ where } \\
 & \mathbf{P_k} = \prod_{i=1}^{m-1}V_{\sigma^k(i)}\prod_{i=1}^{m-1}\Delta_{\sigma^k(i)} = \prod_{i=1,~i\neq k}^m V_i\prod_{i=1,~i\neq k}^m \Delta_i, \\
& N_k = \prod_{i=1}^{m-1} \Delta_{\sigma^k(i)} = \prod_{i=1,~ i\neq k}^{m} \Delta_i, \ \ 1\leq k \leq m.
\end{split}\right.
\end{equation}
Note that both $M_k, N_k \in \mathcal{B}(\W)$ for $1\leq k\leq m$. Also, $\ker(\boldsymbol{\Delta}) = \W$ implies that $N_iN_j = \mathbf{0}$ for all $1\leq i,j \leq m$, $i\neq j$. Product of two lower triangular Toeplitz block matrix is again a lower triangular Toeplitz block matrix. So, for every $\rho\in S_m$, $\displaystyle \prod_{i=1}^m V_{\rho(i)}|_{\X_s}$ is a lower triangular Toeplitz block matrix $\Lambda^{-1}\mathbf{T_L}\left[B_{i-j}^{(m)}\right]\Lambda$ (see (\ref{eq:203}) for notations), where $B_{n}^{(m)}$ ($n\geq 0$) is obtained by the following recurrence relations:
\begin{equation}\label{eq:408}
\left\{\begin{split}
& B^{(1)}_0 = M_{\rho(1)}, \ \  B^{(1)}_1 = N_{\rho(1)} \ \ B^{(1)}_i = \mathbf{0}, \text{ for all } i\geq 2, \\
& B^{(k)}_0 = B^{(k-1)}_0 M_{\rho(k)} \\
& B^{(k)}_i = B^{(k-1)}_i M_{\rho(k)} + B^{(k-1)}_{i-1} N_{\rho(k)}, \ \text{for} 1\leq i \leq k-1, \ k\geq 2, \\
& B^{(k)}_i = \mathbf{0}, \text{ for } \ \  i \geq k \geq 2, \quad \left[\text{ as } N_p N_q = \mathbf{0}, \ \ p \neq q ~\right].
\end{split}\right.
\end{equation}
The relations in (\ref{eq:408}) are obtained by computing the block matrix multiplication from the left, i.e., $(V_{\rho(1)} \dots V_{\rho(m-1)}) V_{\rho(m)}$. If we compute the block matrix multiplication from the right, i.e., $V_{\rho(1)}(V_{\rho(2)}\dotsc V_{\rho(m)})$, then we again reach at the same recurrence relations with different presentation as follows:
\begin{equation*}
\left\{\begin{split}
& A^{(1)}_0 = M_{\rho(m)}; \ \  A^{(1)}_1 = N_{\rho(m)} \ \ A^{(1)}_i = \mathbf{0}, \text{ for all } i\geq 2, \\
& A^{(k)}_0 = A^{(k-1)}_0 M_{\rho(m-k+1)}, \\
& A^{(k)}_i = N_{\rho(m-k+1)}A^{(k-1)}_i  + M_{\rho(m-k+1)}A^{(k-1)}_{i-1}, \ \ \text{for} \ \ 1\leq i \leq k-1, \ k\geq 2, \\
& A^{(k)}_i = \mathbf{0}, \text{ for } \ \  i \geq k \geq 2, \quad \left[\text{ as } N_p N_q = \mathbf{0}, \ \ p \neq q ~\right]
\end{split}\right.
\end{equation*}
Nevertheless, we go ahead with the recurrence relations as in $(\ref{eq:407})$. Since the product $\displaystyle \prod_{i=1}^m V_{\rho(i)}|_{\X_s} = \mathbb{V}|_{\X_s}= \Lambda^{-1} \mathbf{T}_{zI_\W} \Lambda$, we have $B^{(m)}_1 = I_\W$ and $B^{(m)}_i = \mathbf{0}$ for all $i\neq 1$. Thus, we have the following operator equations:
\begin{equation}\label{eq:409}
\left\{\begin{split}
& B^{(m)}_0 = \prod_{i=1}^m M_{\rho(i)} = \mathbf{0},  \\
& B^{(m)}_{i} = B^{(m-1)}_i M_{\rho(m)} + B^{(m-1)}_{i-1} N_{\rho(m)} =  \mathbf{0}, \forall i \neq 1, \\
& \begin{multlined}[t]
B^{(m)}_1= I_\W = \sum_{i=2}^{m-1} \left(M_{\rho(1)}\dotsc M_{\rho(i-1)}N_{\rho(i)}M_{\rho(i+1)}\dotsc M_{\rho(m)}\right)
+ M_{\rho(1)}\dotsc M_{\rho(m-1)}N_{\rho(m)} \\
+ N_{\rho(1)}M_{\rho(2)} \dotsc M_{\rho(m)}, ~~\rho\in S_m.
\end{multlined}
\end{split}\right.
\end{equation}

On account of part-(i) of Lemma \ref{thm:301} and Lemma \ref{lem:412}, we have that $\X_u$ is invariant under each $V_k$ and $V_k|_{\X_u}$ is a unitary (surjective isometry). Consequently, it suffices to consider that $V= \prod_{i=1}^m V_i$ has no unitary part, i.e.,  $V$ is a pure isometry in studying the BCL model for the tuple $(V_1, \dotsc, V_m)$. 

\smallskip

Now, we are in a position to present the main result of this Section, a structure-cum-model theorem for a commuting tuple of Wold isometries that have right-complemented ranges.

\begin{thm}\label{thm:411}
Let $(V_1, \dotsc, V_m)$ be a tuple of commuting Wold isometries with right-complemented ranges on a reflexive, smooth and strictly convex Banach space $\X$ such that $V= \prod_{i=1}^m V_i$ is a pure isometry. Then the following assertions hold:
\begin{enumerate}
\item[(i)] There is a unique closed subspace $\W$ and a tuple $(\mathbf{P_1}, \dotsc, \mathbf{P_m})$ of norm-one projections on $\W$ such that
\begin{equation*} 
\begin{split}
  (V_1, \dotsc, V_m, V) & \equiv  \left(\mathbf{T}_{M_1 + zN_1}, \dotsc, \mathbf{T}_{M_m + zN_m}, \mathbf{T}_{zI_\W} \right) \quad \text{ on } \quad \mathcal{S},\ \ \text{ where } \\
  & \quad  M_k = V_k(I_\W - \mathbf{P_k}), \ \ N_k = \prod_{i=1,~ i\neq k}^m \Delta_i, \ \ 1\leq k \leq m,
\end{split}
\end{equation*}
and $\mathcal{S}$ is as in $(\ref{eq:202})$ for the isometry $V$.

\smallskip

\item[(ii)] The left-inverse of $\mathbf{T}_{M_k + zN_k}$ $(1\leq k\leq m)$ is again a Toeplitz block matrix and is given by $\mathbf{T}_{C_k + \bar{z}D_k}$ where
\begin{align*}
& C_k = \begin{multlined}[t]
 \sum_{i=2}^{m-2} \left(M_{\sigma^k(1)}\dotsc M_{\sigma^k(i-1)}N_{\sigma^k(i)}M_{\sigma^k(i+1)}\dotsc M_{\sigma^k(m-1)}\right) \\
+ M_{\sigma^k(1)}\dotsc M_{\sigma^k (m-2)}N_{\sigma^k (m)} + N_{\sigma^k(1)} M_{\sigma^k (2)}\dotsc M_{\sigma^k (m-1)}, 
\end{multlined} \\
& D_k = \prod_{i=1}^{m-1} M_{\sigma^k(i)}, \quad \text{where} \quad \sigma(i) = i+k, \ \ \text{for} \ \ 1\leq i \leq m-1, \ \sigma(m) =1,
\end{align*}
and summation over empty set is excluded.

\item[(iii)] The operators $\{\Delta_k M_k, M_k \Delta_k:~ 1\leq k \leq m\}$ on $\W$ are bounded projections of norm at most 2 and have norm 1 if and only if 
\[
 \left(\bigoplus_{i=1}^{m-2} \prod_{j=1}^i V_{\sigma^k(j)} \left(\mathcal{L}_{\sigma^k(i+1)} \right) \oplus \mathcal{L}_{\sigma^k(1)} \right) \perp_B \left(\prod_{i=1}^{m-1} V_{\sigma^k(i)}\mathcal{L}_{\sigma^k(m)}\right),~~ 1\leq k\leq m,
\]
where direct sum over empty set is excluded
\end{enumerate}
\end{thm}

\begin{proof}
It follows from (i)-Proposition \ref{prop:411} that $V$ is a Wold isometry with right-complemented range. Let $\X = V(\X) \oplus_\SR \W$ and $V= V|_{\X_u} \oplus V|_{\X_s}$ be the Wold decomposition of $V$. By hypothesis, we have $\X_u = \mathbf{0}$. In addition, by \cite[Lemma 5.1]{JPR}, $V|_{\X_s}$ is a unilateral shift with generating subspace $\W$. Consequently, we have $\X = \oplus_{n=0}^\infty V^n (\W)$ and $V$ is a unilateral shift with right-complemented range $\W$. By Lemma \ref{lem:202}, $\W$ is unique right-complement of $V(\X)$. The proof of (i) follows from (\ref{eq:407}), whereas (ii) is an immediate consequence of block matrix multiplication and the operator relations presented in (\ref{eq:409}). Now, we prove (iii). Note that for $1\leq k \leq m$, $\Delta_k M_k = (I_\W - \mathbf{P_k})$, which is a projection of norm at most 2. We also have, $M_k\Delta_k = V_k (I_\W - \mathbf{P_k}) \Delta_k$. Thus, it is easy to see that $(M_k\Delta_k)^2 = M_k \Delta_k$ for $1\leq k \leq m$. Consequently, $M_k\Delta_k$ is also a bounded projection of norm at most 2 for $1\leq k \leq m$. The range of $(I_\W - \mathbf{P_k})$ is given in (\ref{eq:406}). Consequently, $\Delta_k M_k$ and $M_k \Delta_k$ are norm-one projections if and only if $(I_\W - \mathbf{P_k})(\W) \perp _B \mathbf{P_k}(\W)$. This completes the proof.
\end{proof}

\section[Explicit unitary extension of Wold isometries]{Explicit unitary extension of commuting Wold isometries}\label{Sec:12}

\noindent
Recall that a unilateral shift of some multiplicity on a Hilbert space always extends to a bilateral shift of same multiplicity, e.g., see \cite[Chapter 1, Section 2, Proposition 2.2]{Nagy Foias}. Again, Wold decomposition of a Hilbert space isometry (see \cite{von Neumann, Wold}) splits an isometry into two orthogonal parts of which one is a unitary and the other is a unilateral shift of some multiplicity. So, it follows that every isometry on a Hilbert space can be extended to a unitary on some bigger Hilbert space. This result was further generalized for a family of commuting isometries on a Hilbert space by Brehmer \cite{Brehmer} and It\^o \cite{ITO}, which is stated below.

\begin{thm}[\cite{Brehmer}, \cite{ITO} \& \cite{Nagy Foias}, Chapter 1, Section 9]\label{thm:101}Let $\{V_{\omega}\}_{w\in \Omega}$ be a family of commuting isometries on a Hilbert space $\HS$. Then there is Hilbert space $\mathcal{K}$ containing $\HS$ and a family of unitaries $\{U_{\omega}\}_{\omega\in \Omega}$ on $\mathcal{K}$ such that $U_\omega$ is an extension of $V_\omega$ for every $\omega \in \Omega$.
\end{thm}
One naturally asks if such simultaneous unitary extension is possible for commuting Banach space isometries. Douglas \cite{Douglas} addressed this problem for a semigroup of isometries by exhibiting a unitary representation that extends a given isometric representation on a Banach space. We have already mentioned this in the `Introduction'. Thus, having an appropriate choice of a semigroup, it follows that a commuting $n$-tuple of isometries on a Banach space extends to a commuting $n$-tuple of unitaries on a larger Banach space. This is just an easy observation. However, for the sake of completeness we present an elementary proof to this statement.

\begin{lem}
Let $(V_1, \dotsc, V_n)$ be a commuting tuple of isometries acting on a Banach space $\X$. Then there is a Banach space $\Y\supseteq \X$ and a commuting tuple of unitaries $(U_1, \cdots, U_n)$ on $\Y$ such that $U_i$ extends $V_i$ for $1\leq i \leq n$.
\end{lem}

\begin{proof}
Let $\N_0 = \N \cup\{0\}$. Consider the semigroup $\Sigma:= \N_0 \oplus \cdots \oplus \N_0$ ($n$-times) and the isometric representation $V: \Sigma \to \mathcal{B}(\X)$ defined by
\[
 V(m_1, \dotsc, m_n) = V_1^{m_1} \cdots V_n^{m_n}, \quad (m_1, \dotsc, m_n)\in \Sigma.
\]
Note that $V(e_i)= V_i$, where $e_i\in \Sigma$ is the element whose $i$-th coordinate is $1$ and $0$ elsewhere. It follows from \cite{Douglas} that there exist a Banach space $\Y\supseteq \X$ and a unitary representation $U: \Sigma \to \mathcal{B}(\Y)$ such that $U$ extends $V$; that is $U(m_1, \dotsc, m_n)$ is an extension of $V(m_1, \dotsc, m_n)$ for every $(m_1, \dotsc, m_n)\in \Sigma$. In particular, $U(e_i)$ is an extension of $V(e_i)= V_i$ for $1\leq i \leq n$. The proof is completed by defining the commuting tuple of unitaries $(U_1, \dotsc, U_n)$ by $U_i = U(e_i)$ for $1\leq i \leq n$.
\end{proof}

The aim of this Section is to construct explicitly such a unitary extension for a certain class of commuting finite isometric tuples. We first show the construction for $n=1,2$ cases without jumping directly to an arbitrary $n$-variable case as we believe that it will ease computational complications and a reader's intuition will be built this way. Also, the construction in the general case is a straight-forward analogue of the $n=2$ case, which is why we refrain from giving an explicit construction in $n$-variable setting. Notably, the passage from $n=1$ to $n=2$ is built on additional assumptions on the concerned isometries and the Banach space they act on. Our algorithm consists of the following steps:
\begin{itemize}
\item for an isometry $V$ with right-complemented range (not necessarily a Wold isometry) on a Banach space $\X$, we find an invertible operator $U$ on a larger Banach space that extends $V$;

\smallskip

\item we find a necessary and sufficient condition on $V$ and $\X$ such that $U$ becomes a unitary extension of $V$;

\smallskip

\item we explain in detail why the same method does not work in general for a pair of commuting isometries with right-complemented ranges on any Banach space and why we need to choose Wold isometries acting on an appropriate Banach space;

\smallskip

\item we show an explicit construction of a commuting pair of invertible operators $(U_1,U_2)$ that extends a given pair of commuting Wold isometries with right-complemented ranges $(V_1,V_2)$ acting on a reflexive, smooth and strictly convex Banach space;

\smallskip

\item a characterization analogous to the one-variable case is found for such pair of isometries $(V_1,V_2)$ and the underlying Banach space so that the concerned extension $(U_1,U_2)$ becomes a unitary extension of $(V_1,V_2)$;

\smallskip

\item it will be evident that the same mechanism (as in the construction in n = 2 case) can be applicable to a general $n$-variable case when $n\geq 2$ and consequently a generalized version of unitary extension result is stated in the form of Theorem \ref{thm:503}.

\end{itemize}
As is mentioned above, a necessary and sufficient condition is given on the isometric tuple $(V_1, \dots , V_n)$ and on the underlying Banach space $\X$ such that the explicitly constructed invertible extension becomes a unitary extension of $(V_1, \dots , V_n)$. Indeed, this characterization is obtained in terms of $2$-direct sum complement of the range of each $V_i$ (see the statement of Theorem \ref{thm:503}), which holds naturally if each $V_i$ acts on a Hilbert space. The other assumptions on the isometries (i.e., each of them is a Wold isometry with right-complemented range) and on $\X$ (i.e., reflexivity, smoothness and strict convexity) are automatic in Hilbert space setting. Hence, our unitary extension generalizes the existing theorem on unitary extension of a commuting tuple of Hilbert space isometries.

\smallskip

In \cite{JPR} authors introduce $\sigma$-shift between Banach spaces in order to generalize Hilbert space bilateral shift. A $\sigma$-shift is a Banach space unitary by definition and it becomes a Hilbert space bilateral shift when acts on a Hilbert spaces as shown in \cite[Proposition 6.6]{JPR}. In \cite[Theorem 6.7]{JPR}, they prove by an explicit construction that every unilateral shift (see Definition \ref{defn:Banach-Unilateral}) on a Banach space extends to a $\sigma$-shift and thus generalize the same result in Hilbert space theory. More generally, in \cite[Theorem 6.9]{JPR}, they show that every Wold isometry on a Banach space extends to a Banach space unitary. However, these unitary extensions act between two different Banach spaces. In this Section, we choose all components of a commuting Banach space isometric tuple to act on one Banach space. Thus, the previous approach (as in \cite[Theorem 6.7, Theorem 6.9]{JPR}) does not work here. Let us unfold our algorithm with the following lemma.

\begin{lem}\label{Lem:14-07-26::19:52}
Let $V$ be an isometry on a Banach space $\X$ with right-complemented range. Let $\mathcal{L}$ be a right-complement of $V(\X)$ in $\X$, and $\Y$ be the Banach space given by $\Y=\ell_2(\mathcal{L})\bigoplus_2 \X$. Then the operator $U: \Y \to \Y$ defined by
\begin{equation}\label{Eq:14-07-26::20:02}
  U(\{l_n\}, x) = (\{l_{n+1}\}, Vx + l_1), \quad x\in \X,~\{l_n\} \in \ell_2(\mathcal{L}), ~l_n \in \mathcal{L}\text{ for all } n\in \mathbb{N},
\end{equation}
is invertible and extends $V$.
\end{lem}

\begin{proof}
Let $(\{l_n\}, x) \in \Y$ be arbitrary. Then we have
\begin{equation}\label{eq:501}
   \| U(\{l_n\}, x) \| = \left( \sum_{n=2}^\infty \|l_n \|^2  + \| Vx + l_1 \|^2 \right)^{\frac{1}{2}} \leq \left( \sum_{n=1}^\infty \|l_n \|^2 \right)^{\frac{1}{2}}  + \| Vx \| \leq 2 \|(\{l_n\}, x) \|.
\end{equation}
So, $U$ is a well defined bounded linear map on $\Y$. The injectivity of $U$ follows from its definition. On the other hand, for every $\mathbf{y}= (\{r_n\}, y)\in \Y,$ there is $x\in \X$ and $l\in \mathcal{L}$ such that $y = Vx + l.$ Thus,
\[
 U(\{r_n'\}, x) = \mathbf{y}, \quad \text{ where } \quad r_n' = \begin{cases}
                                             l & \text{ if } n=1 \\
                                             r_{n-1} & \text{ if } n \geq 2
                                           \end{cases}\quad \text{ and } \quad \{r'_n\}\in \ell_2(\mathcal{L}).
\]
This shows that $U$ is surjective. Therefore, by open mapping theorem, $U^{-1}$ is also a bounded linear map. The natural inclusion map $W: \X \to \Y$ given by $x \mapsto (\mathbf{0}, x)$ is an isometric embedding of $\X$ into $\Y$. The operator $U$ extends $V$ because $UW(x) = WV(x)$ for all $x\in \X$. This completes the proof.
\end{proof}

Note that the operator $U$ defined in Lemma \ref{Lem:14-07-26::19:52} by Equation-\eqref{Eq:14-07-26::20:02} depends on the right-complement $\mathcal{L}$ of $V(\X)$, which is not necessarily unique on an arbitrary Banach space $\X$. Moreover, this invertible extension is not a unitary in general. As mentioned before, we now find a necessary and sufficient condition such that the invertible extension $U$ is a unitary.

\begin{thm}\label{thm:14-07-26::20:15}
The operator $U$ as in the statement of $Lemma \ \ref{Lem:14-07-26::19:52}$ is a unitary if and only if $\X = V(\X)\bigoplus_2 \mathcal{L}$ for some closed subspace $\mathcal{L}$ of $\X$.
\end{thm}

\begin{proof}

First we prove the necessity. Suppose $U$ is a unitary. Let $l\in \mathcal{L}$ and $x\in \X$ be arbitrary. Set $l_1= l$ and $l_n = \mathbf{0}$ for all $n\geq 2$. Then $\|U(\{l_n\},x)\| = \|(\{l_n\},x)\|$ implies that $\|Vx+l\|^2 =\|Vx\|^2 + \|l\|^2$. Consequently, we have $\X = V(\X) \oplus_2 \mathcal{L}$.

\smallskip

Now, we prove the sufficiency. Let $\widehat{\mathcal{L}}$ be a subspace such that $\X = V(\X) \oplus_2 \widehat{\mathcal{L}}$. Consider $\mathcal{L}= \widehat{\mathcal{L}}$ and define the operator $U$ accordingly as in \eqref{Eq:14-07-26::20:02}. Then it follows from (\ref{eq:501}) that for all $(\{l_n\},x)\in \Y$ we have
\[
\| U(\{l_n\}, x) \|^2 = \sum_{n=2}^\infty \|l_n \|^2  + \| Vx + l_1 \|^2 = \sum_{n=1}^\infty \|l_n \|^2  + \|x\|^2 = \|(\{l_n\}, x) \|^2.
\]
Therefore, $U$ is an isometry. Also, it follows from Lemma \ref{Lem:14-07-26::19:52} that $U$ is invertbile. Hence $U$ is a unitary. This completes the proof.
\end{proof}

The assumptions in Lemma \ref{Lem:14-07-26::19:52} and Theorem \ref{thm:14-07-26::20:15} are automatically satisfied when $\X$ is a Hilbert space. Now, suppose $n \geq 2$ and $(V_1, \dotsc, V_n)$ is a commuting tuple of isometries on a Banach space $\X$ such that the range of each $V_i$ is right-complemented in $\X$. To construct a Banach space $\Y$ (as in Lemma \ref{Lem:14-07-26::19:52}) and a commuting tuple $(U_1, \dotsc, U_n)$ of invertible operators on $\Y$, the first step is to determine a right-complement of the isometry $V=\prod_{i=1}^n V_i$. The proof of  part (ii) in Proposition \ref{prop:411} shows that a right-complement of the isometry $V=\prod_{i=1}^n V_i$ can be obtained in the same way under only smoothness assumption on $\X$. To see this, let us suppose $n=2$ and let $\X$ be smooth. Then, the proof of part (i) in Lemma \ref{lem:411} shows that both $V_1(\mathcal{L}_2) \oplus_\SR \mathcal{L}_1$ and $ V_2(\mathcal{L}_1) \oplus_\SR \mathcal{L}_2$ are right-complements of $V_1V_2(\X)$, where $\mathcal{L}_i$ is a right-complement of $V_i(\X)$ for $i=1,2$. Moreover, the assumptions of reflexivity, smoothness and strict convexity on the underlying Banach space give a unique right-complement of the isometry $V= \prod_{i=1}^n V_i$ as shown in Proposition \ref{prop:411}. Thus, one can suitably define the space $\Y$. Once the space $\Y$ is fixed, the next step is to construct a commuting tuple $(U_1, \dotsc, U_n)$ of invertible operators on $\Y$ such that $U_i$ is an extension of $V_i$ for each $i$. At this point, we bring into the picture the additional assumption of each $V_i$ being a Wold isometry. Note that a Wold isometry with right-complemented range acting on a reflexive, smooth and strictly convex Bnach space possesses a unique left-inverse (by Lemma \ref{lem:202}). We use this property to show that the operators $U_i$, $1\leq i \leq n$, are well defined and commutative. Needless to mention that the adjoint of a Hilbert space isometry $\widehat{V}$ plays the role of left-inverse and thus our theory generalizes the corresponding Hilbert space result. 

\smallskip

Let $V$ be a Wold isometry with right-complemented range on a Banach space $\X$ with a Wold decomposition $V= V|_{\X_u} \oplus V|_{\X_s}$ as in \eqref{Wold decomposition Banach spaces}. Then by \cite[Lemma 5.1]{JPR}, $V|_{\X_s}$ is a unilateral shift with generating subspace $\mathcal{L}$. Consider the space $\mathbb{G}= \X_u\oplus \mathcal{S}$, where $\mathcal{S}$ is as in (\ref{eq:201}). Then $\X$ is isometrically isomorphic with $\mathbb{G}$ via the unitary $\Lambda_v: \X \to \mathbb{G}$ defined by 
\begin{equation}\label{eq:502}
\Lambda_v\left(x_u + \sum_{n=0}^\infty V^n l_n\right) = \left(x_u, \{l_n\}_{n\in \N_0}\right).
\end{equation}
The isometry $V$ can be realised as a Wold isometry with right-complemented range $\widehat{V}$ on $\mathbb{G}$, where $\widehat{V} = \Lambda_v V \Lambda_v^{-1}$. Let us follow this and as mentioned before, let us consider the $n=2$ case first.

\begin{thm}\label{thm:502}
Let $(V_1,V_2)$ be a pair of commuting Wold isometries with right-complemented ranges acting on a reflexive, smooth and strictly convex Banach space $\X$. Then there is a pair $(U_1, U_2)$ of commuting invertible operators acting on a Banach space $\Y$ containing $\X$ and $U_1^{n_1}U_2^{n_2}$ is an extension of $V_1^{n_1}V_2^{n_2}$ for all $n_1$, $n_2 \geq 0$. Moreover, $(U_1,U_2)$ is a pair of commuting unitaries if and only if $ \X= V_i\X \bigoplus_2 \mathcal{L}_i$ for some closed subspaces $\mathcal{L}_i$ for $i=1,2$. Therefore, $U_1, U_2$ are unitaries when $\X$ is a Hilbert space.
\end{thm}

\begin{proof}
Let $\mathcal{L}_i$ be the right-complement of $V_i(\X)$ and $\Delta_i$ be the associated left-inverse of $V_i$ for $i=1,2$. Also, let $V= V_1 V_2$. Then $V$ is a Wold isometry with right-complemented range by (i)-Lemma \ref{lem:411} and the right-complement $\W$ of $V(\X)$ is given by $\W = V_1\mathcal{L}_2 \oplus_{\SR} \mathcal{L}_1 = V_2\mathcal{L}_1 \oplus_\SR \mathcal{L}_2$. Let $V= V|_{\X_u} \oplus V|_{\X_s}$ be the Wold decomposition of $V$. The space $\mathbb{G} = \X_u \bigoplus \mathcal{S}$ is isometrically isomorphic with $\X$ via $\Lambda_v$ (see (\ref{eq:502})). Consider the Banach space $\Y = \ell_2(\W) \oplus_2\mathbb{G}$. Let $\phi: \mathbb{G} \to \Y$ be the natural inclusion map $g\mapsto (\mathbf{0}, g)$. Then $\phi\circ \Lambda_v: \X \to \Y$ is an isometric embedding of $\X$ into $\Y$. A general element in $\Y$ is of the form $(\{w_n\}_{n\in \Z_-}, (x_u, \{w_n\}_{n\in \N_0}))$. To avoid notation complication, we write $(\{w_n\}_{n\in \Z_-}, (x_u, \{w_n\}_{n\in \N_0}))$ as $(\{w_n\}_{n\in \Z}, x_u)$, where it is understood that the sequence $\{w_n\}_{n\in \Z_-}$ is a member of $\ell_2(\W)$ and $\{w_n\}_{n\in \N_0}$ is a member of $\mathcal{S}$. Let $(\{w_n\}_{n\in \Z}, x_u) \in \Y$ be arbitrary. Then there are unique sequences $\{l_n(1)\}_{n\in \Z}, \{r_n(1)\}_{n\in \Z}$ in $\mathcal{L}_1$ and $\{l_n(2)\}_{n\in \Z}, \{r_n(2)\}_{n\in \Z}$ in $\mathcal{L}_2$ such that 
\[
V_2 l_n(1) + l_n(2) = w_n = V_1r_n(2) + r_n(1), \quad n\in \mathbb{Z}.
\]
Now, define $U_1,~ U_2: \Y \to \Y$ by
\begin{align*}
    & U_1 (\{w_n\}_{n\in \Z}, x_u) = (\{w_n'\}_{n\in \Z}, V_1x_u), \quad w_n' = V_1 l_n(2) + l_{n-1}(1), \quad n \in \mathbb{Z}, \\
    & U_2 (\{w_n\}_{n\in \Z}, x_u) = (\{w_n''\}_{n\in \Z}, V_2x_u), \quad w_n'' = V_2 r_n(1) + r_{n-1}(2), \ n \in \mathbb{Z}.
\end{align*}
We show in a few steps below that $U_1, U_2$ are well-defined, bounded linear operator, commuting, invertible and $U_1^{m_1}U_2^{m_2}$ extends $V_1^{m_1}V_2^{m_2}$ for all $m_1, m_2 \geq 0$.

\smallskip

\noindent $\bullet$ We show that $U_i$ is well-defined and linear for $i=1,2$. Let $P_i = (I_\W - Q_{i})$, where $Q_i$ ($i=1,2$) is the norm-one projection on $\W$ defined in (iv)-Lemma \ref{lem:411}. Then we have
\begin{equation}\label{eq:503}
  U_i (\{w_n\}_{n\in \Z}, x_u) = \left(\{V_iP_jw_n+ \Delta_jw_{n-1}\}_{n\in \Z}, V_ix_u \right), \quad i,j\in \{1,2\}, \ i\neq j.
\end{equation}
It follows from (\ref{eq:503}) that each $U_i$ is linear. We only prove that $U_1$ is well-defined, the proof for $U_2$ is similar. Note that $\{V_1P_2w_n + \Delta_2w_{n-1}\}_{n\in \Z_-} \in \ell_2(\W)$ because
\begin{align*}
  \left(\sum_{n=-\infty}^{-1} \|V_1P_2w_n + \Delta_2w_{n-1}\|^2 \right)^{\frac{1}{2}} & \leq 
  \left( \sum_{n=-\infty}^{-1} \|V_1P_2w_n\|^2 \right)^{\frac{1}{2}} + \left( \sum_{n=-\infty}^{-1} \|\Delta_2w_{n-1}\|^2 \right)^{\frac{1}{2}}\\
& \leq 2 \left(\sum_{n=-\infty}^{-1} \|w_n\|^2 \right)^{\frac{1}{2}} + \left(\sum_{n=-\infty}^{-1} \|w_{n-1}\|^2 \right)^{\frac{1}{2}} \\
& \leq 3 \left\|\{w_n\}_{n\in \Z_-}\right\| < \infty,
\end{align*}
where the second inequality follows from $\|P_2\|\leq 2$ and $\|\Delta_2\| =1$. It remains to show that the series $\displaystyle \sum_{n=0}^\infty V^n (V_1P_2(w_n) + \Delta_2(w_{n-1}))$ converges in $\X_s$. Let $\{s_n\}$ be the associated partial sum sequence. Then for every $ n\geq 1$, $s_n$ can be expressed as
\begin{align}\label{eq:504}
  s_n & = \begin{multlined}[t] 
  l_{-1}(1) + V_1 \left(V_2l_0(1) + l_0(2)\right) + V_1 V\left( V_2 l_0(1) + l_0(1)\right) + \cdots + \\
  + V_1 V^{n-1}\left(V_2l_n(1) + l_n(2)\right) + V_1 V^n l_n(2) 
  \end{multlined}\nonumber\\
     & = l_{-1}(1) + V_1 \left(\sum_{k=0}^n V^k (V_2 l_k(1) + l_k(2)) \right) + V_1 V^n l_n(2) \nonumber \\
     & = l_{-1}(1) + V_1 \left(\sum_{k=0}^n V^k w_k \right) + V_1 V^n l_n(2).
\end{align}
The convergence of $\displaystyle \sum_{n=0}^\infty V^n w_n $ implies that $\displaystyle \lim_{n\to \infty} V^n w_n = \mathbf{0}$. So, $\displaystyle \lim_{n\to \infty} \|V^n w_n\| = \lim_{n\to \infty}\|w_n\| = 0$. Thus, it follows from (\ref{eq:504}) that $\{s_n\}$ is a convergent sequence, since $\displaystyle \lim_{n\to \infty} V_1 V^n l_n(2) = \mathbf{0}$. Indeed,
\[
   \| V_1 V^n l_n(2) \| = \| V_1 V^n P_2 w_n \|\leq 2 \| w_n \| \longrightarrow 0 \ \text{ as } n \to \infty.
\]
Moreover, we have
\begin{equation}\label{eq:505}
  \sum_{n=0}^\infty V^n (V_1 P_2 w_n + \Delta_2 w_{n-1}) = \lim_{n\to \infty}s_n = V_1\left(\sum_{n=0}^\infty V^n w_n \right) + l_{-1}(1).
\end{equation}

\noindent $\bullet$ Now, we show that $U_i$ is bounded for $i=1,2$ with an estimate of an upper bound for norm of each $U_i$. Let $(\{w_n\}_{n\in \Z}, x_u) \in \Y$ be arbitrary. Then, we have
\begin{align*}
\left\| U_1\left(\{w_n\}_{n\in \Z}, x_u\right) \right\| & = 
\left[ \sum_{n=-\infty}^{-1} \left\| V_1 P_2 w_n + \Delta_2 w_{n-1} \right\|^2
+ \left \| \left(\{V_1P_2 w_n + \Delta_2 w_{n-1}\}_{n\in \N_0}, V_1 x_u \right) \right\|^2 \right]^{1\slash 2} \\
 & = \left[ \sum_{n=-\infty}^{-1} \left\| V_1 P_2 w_n + \Delta_2 w_{n-1} \right\|^2 + \left\| V_1 x_u + \sum_{n=0}^\infty V^n (V_1P_2 w_n + \Delta_2 w_{n-1}) \right\|^2 \right]^{1\slash 2} \\
  & = \left[ \sum_{n=-\infty}^{-1} \| V_1 P_2 w_n + \Delta_2 w_{n-1} \|^2 + \left\| V_1 \left( x_u + \sum_{n=0}^\infty V^n w_n \right) + \Delta_2w_{-1}\right\|^2 \right]^{\frac{1}{2}} ~\left[ \text{ by } \eqref{eq:505} \right] \\
  & \leq \begin{multlined}[t]
  \left[ \sum_{n=-\infty}^{-1} \left\|V_1 P_2 w_n \right\|^2 + \left\| V_1\left( x_u +  \sum_{n=0}^\infty V^n w_n \right) \right\|^2 \right]^{1\slash 2} \\
  + \left[ \sum_{n=-\infty}^{-1} \| \Delta_2 w_{n-1} \|^2 + \| \Delta_2 w_{-1} \|^2 \right]^{1\slash 2} 
  \end{multlined}\\
 & \leq \left[ \sum_{n=-\infty}^{-1} 4 \| w_n \|^2 + \left\| x_u + \sum_{n=0}^\infty V^n w_n \right\|_\X^2 \right]^{\frac{1}{2}} + \left[ \sum_{n=-\infty}^{-1} \| w_n\|^2 \right]^{\frac{1}{2}} \\
 & \leq 3 ~ \left\| \left(\{w_n\}_{n\in \Z}, x_u\right) \right\|,
\end{align*}
where the second last inequality follows from the facts that $\|P_2\|\leq 2$ and $\|\Delta_2\|=1$. Similarly, $\left\| U_2(\{w_n\}_{n\in \Z}, x_u) \right\| \leq 3 \left\|(\{w_n\}_{n\in \Z}, x_u)\right\|$ for every $(\{w_n\}_{n\in \Z}, x_u) \in \Y$. Consequently, we have $\|U_i\| \leq 3$ for $i=1,2$.

\smallskip

\noindent $\bullet$ We show that each $U_i$ is boundedly invertible. By open mapping theorem it suffices to show that $U_i$ is bijective. The injectivity of $U_i$ ($i=1,2$) follows from its definition. We show the surjectivity of $U_1$. The surjectivity of $U_2$ follows in a similar manner. Let $ \mathbf{y}= (\{w_n\}_{n \in\Z}, x_u)\in \Y$ be arbitrary. Then there exist sequences $\{l_n(i)\}_{n\in \Z}$ in $\mathcal{L}_i$ ($i=1,2$) such that $w_n = V_1 l_n(2) + l_n(1)$ for all $n \in \mathbb{Z}$. Also, since $V_1|_{\X_u}$ is unitary, there is $x_{u_1}\in \X_u$ such that $V_1(x_{u_1}) = x_u$. Now, consider the element $(\{w'_n\}_{n\in \Z}, x_{u_1})$ defined by $w'_n = V_2l_{n+1}(1) + l_n(2)$ for all $n\in \mathbb{Z}$. It is easy to see that $U_1(\{w'_n\}_{n\in \Z},x_{u_1}) = \mathbf{y}$. It remains to show that $\{w_n'\}_{n\in \Z_-}\in \ell_2(\W)$ and $\{w_n'\}_{n\in \N_0}\in \mathcal{S}.$ Indeed,
\begin{align*}
    \left(\sum_{n=-\infty}^{-1} \| w'_n\|^2 \right)^{\frac{1}{2}} & = \left( \sum_{n=-\infty}^{-1} \| V_2 P_1 w_{n+1} + \Delta_1 w_n \|^2 \right)^{\frac{1}{2}} \\
    & \leq \left( \sum_{n=-\infty}^{-1} 4 \| w_{n+1}\|^2 \right)^{\frac{1}{2}} + \left( \sum_{n=-\infty}^{-1} \| w_n \|^2 \right)^{\frac{1}{2}} < \infty,
\end{align*}
and the convergence of $\displaystyle\sum_{n=0}^\infty V^n w_n'$ is identical to the proof of  (\ref{eq:505}).

\smallskip

\noindent $\bullet$ To show that $U_1,U_2$ commute. By (\ref{eq:503}), we have
\[
  U_1 U_2(\{w_n\}_{n\in \Z}, x_u) = U_1 (\{V_2 P_1 w_n + \Delta_1 w_{n-1}\}_{n\in \Z}, V_2 x_u) = (\{z_n\}_{n\in \Z}, Vx_u), 
\]
where $\{z_n\}_{n\in \Z}$ is given by
\begin{align*}
  z_n & = V_1 P_2 (V_2 P_1 w_n + \Delta_1 w_{n-1}) + \Delta_2(V_2 P_1 w_{n-1} + \Delta_1 w_{n-2}), \ n\in \mathbb{Z} \\
      & = V_1 P_2V_2 r_n(1) + V_1 P_2 r_{n-1}(2) + P_1w_{n-1} + \Delta_2\Delta_1 w_{n-2} \\
      & = V_1 r_{n-1}(2) + r_{n-1}(1) + \Delta_2 \Delta_1 w_{n-2} \\
      & = w_{n-1}, \quad n \in \mathbb{Z}.
\end{align*}
The third equality follows from the fact that $\ker(P_2)= V_2\mathcal{L}_1,  Ran(P_2) = \mathcal{L}_2$ and the last equality holds as $\Delta_2\Delta_1(\W) = \{\mathbf{0}\}$. Similarly, by $\Delta_1\Delta_2 (\W) =\{\mathbf{0}\}$ we have
\begin{equation}\label{eq:506}
U_2 U_1(\{w_n\}_{n\in \Z}, x_u) = (\{w_{n-1}\}_{n\in \Z}, Vx_u) = U_1 U_2 (\{w_n\}_{n\in \Z}, x_u).
\end{equation}
Consequently, $(U_1,U_2)$ is a commuting pair of invertible operators acting on $\Y$. 

\medskip

\noindent $\bullet$ We show that $U_1^{m_1}U_2^{m_2}$ extends $V_1^{m_1}V_2^{m_2}$ for all $m_1, m_2\geq 0$. A first step towards this will be proving $U_i$ extends $V_i$ for $i=1,2$. For any $x\in \X$, there are unique $x_u\in \X_u$ and $x_s\in \X_s$ such that $x= x_u + x_s$. So, there is a unique sequence $\{w_n\}_{n\in \N_0}$ in $\W$ such that $x_s = \sum_{n=0}^\infty V^n w_n$. For every $n\geq 0$, there are unique $\{l_n(i), r_n(i)\} \subseteq \mathcal{L}_i$ ($i=1,2$) such that $w_n = V_2 l_n(1) + l_n(2) = V_1r_n(2) + r_n(1) $ for all $n \geq 0$. So, we have
\begin{align}
   V_1(x_s) & = V_1l_0(2) + \sum_{n=1}^\infty V^n (V_1l_n(2) + l_{n-1}(1)) = V_1l_0(2) + \sum_{n=1}^\infty V^n (V_1P_2w_n + \Delta_2w_{n-1}) \label{eq:507}
\end{align}
and 
\begin{align}
   V_2(x_s) & = V_2 r_0(1) + \sum_{n=1}^\infty V^n (V_2r_n(1) + r_{n-1}(2)) = V_2 r_0(1) + \sum_{n=1}^\infty V^n (V_2P_1w_n + \Delta_1w_{n-1}), \label{eq:508}
\end{align}
where the proof of convergence of both the series $\sum_{n=1}^\infty V^n (V_1l_n(2) + l_{n-1}(1))$ and $\sum_{n=0}^\infty V^n (V_2r_n(1) + r_{n-1}(2))$ are similar to the proof of (\ref{eq:505}). Denote by $\psi= \phi\circ \Lambda_v$ to be the isometric embedding of $\X$ into $\Y$. Let $\psi(x)=(\{\tau_n\}_{n\in \Z}, y) $. Then $y= x_u$, $\tau_n = w_n$ for all $n\geq 0$ and $\tau_n = \mathbf{0}$ for all $n\leq -1$. Therefore, the definition of $U_1, U_2$ and (\ref{eq:507}), (\ref{eq:508}) implies that
\[
   U_1 \psi(x) = \psi V_1(x) \quad \text{ and } \quad U_2 \psi(x) = \psi V_2(x).
\]
Hence, $U_i$ extends $V_i$ for $i=1,2$. Now, we show that $U_1^{m_1}U_2^{m_2}$ extends $V_1^{m_1}V_2^{m_2}$ for all $m_1, m_2\geq 0$. Again, (\ref{eq:507}) implies that
\begin{align}
  V_2V_1(x_s) & = V_2V_1l_0(2)  + VV_2l_0(1) + \sum_{n=2}^\infty V^n (V_2l_{n-1}(1) + l_{n-1}(2)) \nonumber \\
  & = \sum_{n=1}^\infty V^n (V_2l_{n-1}(1)) + l_{n-1}(2)) = \sum_{n=1}^\infty V^n w_{n-1} \label{eq:509}.
\end{align}
Therefore, by (\ref{eq:506}) and (\ref{eq:509}) we have 
\[
U_1U_2(\psi(x)) = U_2U_1(\psi(x)) = \psi(V_1V_2(x)) = \psi(V_2V_1(x)).
\]
So, for all $m_1, m_2 \geq 1$ we have 
\begin{align}
U_1^{m_1}U_2^{m_2}(\psi(x)) & = U_1^{m_1-1}U_2^{m_2-1}U_1U_2(\psi(x)) \nonumber \\
& = U_1^{m_1-1}U_2^{m_2-1}(\psi(V_1V_2(x)) = \cdots = \psi(V_1^{m_1}V_2^{m_2}(x)).
\end{align}
Consequently, $U_1^{m_1}U_2^{m_2}$ is an extension of $V_1^{m_1}V_2^{m_2}$ for all $m_1, m_2 \geq 0$.

\medskip

\noindent \textbf{Proof of the second part.} We show the sufficiency first. Suppose $\X = V_i(\X)\oplus_2 \mathcal{L}_i$ for $i=1,2$. Then the right-complement $\W$ of $V(\X)$ is given by $V_1\mathcal{L}_2 \oplus_2 \mathcal{L}_1 = \W = V_2\mathcal{L}_1 \oplus_2 \mathcal{L}_2$. We show that $U_1$ is unitary. The proof for $U_2$ is similar. Let $(\{w_n\}_{n\in \Z}, x_u)\in \Y$ be arbitrary such that $w_n= V_2 l_n(1) + l_n(2)$, where $l_n(i)\in \mathcal{L}_i$ for all $n\in \mathbb{Z}$, $i=1,2$. Now, it follows from (\ref{eq:505}) that
\[
 \left\|U_1\left(\{w_n\}_{n\in \Z}, x_u\right)\right\|^2  = \sum_{n=-\infty}^{-1}\|w_n'\|^2 + \left\|V_1(x_u) + l_{-1}(1) + V_1\left(\sum_{n=0}^\infty V^n w_n\right)\right\|^2,
\]
where $w_n' = V_1 l_n(2) + l_{n-1}(1)$. Thus, we have
\begin{align*}
 \left\|U_1\left(\{w_n\}_{n\in \Z}, x_u\right)\right\|^2 & = \sum_{n=-\infty}^{-1} (\|l_n(2)\|^2 + \|l_{n-1}(1)\|^2 ) + \|l_{-1}(1)\|^2 + \left\|V_1\left(x_u + \sum_{n=0}^\infty V^n w_n\right)\right\|^2 \\
 & = \left\|\{w_n\}_{n\in \Z_-}\right\|^2 + \left\|x_u + \sum_{n=0}^\infty V^n w_n \right\|^2 \\
 & = \left\|(\{w_n\}_{n\in \Z}, x_u )\right\|.
\end{align*}
Therefore, $U_1$ is an isometry. Hence $(U_1, U_2)$ is a pair of commuting unitaries.

\smallskip

Now, we prove the converse, i.e., the necessity. Suppose $(U_1, U_2)$ is a pair of commuting unitaries. Then $U_1U_2$ is a unitary extension of $V$. Therefore, by Theorem \ref{thm:14-07-26::20:15}, we have $\X = V(\X) \oplus_2 \W$. We show that $\X= V_1(\X) \oplus_2 \mathcal{L}_1$. The proof for $V_2$ is similar. Let $x\in \X$ and $l_1\in \mathcal{L}_1$ be arbitrary. Then there are unique $x_u\in \X_u$ and $w\in \W$ such that $x= x_u + w$. Since $V|_{\X_u}$ is unitary, $V_1x_u \in \X_u$. So, there are unique vectors $r_i\in \mathcal{L}_i$ ($i=1,2$) such that $w= V_2r_1 + r_2$. Therefore,
\begin{equation}\label{eq:510}
  \|V_1x + l_1\| = \|V_1(x_u + V_2r_1 + r_2) + l_1\| = \|V_1x_u + V_1V_2 r_1 + V_1r_2 + l_1 \|.
\end{equation}
Define $w_0 = V_2r_1 + r_2$, $w_{-1}= V_2l_1$ and $w_n = \mathbf{0}$ for all $n\in \mathbb{Z}\setminus\{-1,0\}$. Then $\left(\{w_n\}_{n\in \Z}, x_u\right)\in \Y$. Since $U_1$ is an isometry, we have 
\[
 \left\| U_1(\{w_n\}_{n\in \Z}, x_u) \right\|^2 = \left\|(\{w_n\}_{n\in \Z}, x_u)\right\|^2,
\]
which implies that
\[
 \| V_1x_u + l_1 + V_1r_2 + V_1 V_2 r_1 \|^2 = \|V_2l_1\|^2 + \|x_u + V_2r_1 + r_2\|^2 = \|l_1\|^2 + \|V_1x\|^2,
\]
where the last equality follows from the facts that $V_1, V_2$ are isometries and $x= x_u + w$. Therefore, by (\ref{eq:510}) we have $\| V_1x + l \|^2 = \|V_1x\|^2 + \|l\|^2$. Consequently, $\X = V_1(\X) \oplus_2 \mathcal{L}_1$. This completes the proof.
\end{proof}

Now, let us consider the general case. Let $(V_1, \dotsc, V_n)$ be a commuting tuple of Wold isometries with right-complemented ranges acting on a reflexive, smooth, and strictly convex Banach space $\X$. Let $\mathcal{L}_i$ denote the right-complement of $V_i(\X)$ in $\X$ for each $i$. We proved in Proposition \ref{prop:411} that $V=\prod_{i=1}^n V_i$ is a Wold isometry with right-complemented range, say $\W$. In part (ii) of Proposition \ref{prop:411}, we also found the right-complement $\W$ and expressed it corresponding to each permutation $\sigma \in S_n$. A closer look at the construction in the $n=2$ case reveals that the representation of $\W$ corresponding to the permutation $(2~1)$ is used to construct $U_1$, and the identity permutation $(1~2)$ is considered while constructing $U_2$. In the general case, we use the representations of $\W$ in a similar fashion as we have done in Section \ref{Sec:11} to establish a BCL model. The remaining calculations are analogous to those of the $n=2$ case. To avoid unnecessary complications, we refrain from showing such notoriously elaborative and complicated computations of the general case. In conclusion, we present the main result of this Section.

\begin{thm}\label{thm:503}
Let $n\geq 2$ and let $(V_1, \dotsc, V_n)$ be a tuple of Wold isometries with right-complemented ranges acting on a reflexive, smooth and strictly convex Banach space $\X$. Then there is a Banach space $\Y\supseteq\X$ and a tuple $(U_1, \dotsc, U_n)$ of commuting invertible operators on $\Y$ such that $U_1^{m_1} \dotsc U_n^{m_n}$ is an extension of $V_1^{m_1} \dotsc V_n^{m_n}$ for all $m_i\geq 0$, $1\leq i \leq n$. Moreover, $(U_1,\dotsc, U_n)$ is a tuple of commuting isometries if and only if $\X = V_i (\X) \oplus_2 \mathcal{L}_i$ for some closed subspaces $\mathcal{L}_i$ for $1\leq i \leq n$.
\end{thm}

As mentioned before, Theorem \ref{thm:503} generalizes the celebrated theorem due to Brehmer \cite{Brehmer} and It\^o \cite{ITO} on simultaneous unitary extension of commuting Hilbert space isometries. Hence, we state the concerned result in Hilbert space setting as a corollary of Theorem \ref{thm:503}.

\begin{cor}
Let $n\geq 2$ and let $(V_1, \dotsc, V_n)$ be a tuple of isometries acting on a Hilbert space $\HS$. Then there is a Hilbert space $\mathcal{K}\supseteq\HS$ and a tuple $(U_1,\dotsc, U_n)$ of commuting unitaries on $\mathcal{K}$ such that $U_1^{m_1}\dotsc U_n^{m_n}$ is an extension of $V_1^{m_1}\dotsc V_n^{m_n}$ for all $m_i\geq 0$, $1\leq i \leq n$.
\end{cor}

We end this paper with the positive outcome that our construction of unitary extension generalizes that of Hilbert space isometries. It is possible mainly because the concerned unitary tuple acts on a Banach space which is an orthogonal $2$-sum of the given Banach space and some other Banach space. However, it is to be noted that Douglas's result on unitary extension \cite{Douglas} is valid for any commuting tuple of Banach space isometries even though it is existential, whereas we work only with Wold isometries having right-complemented ranges and acting on a reasonably nice Banach space. Also, our extension is first invertible and then is conditionally unitary and acts on a Banach space of specific kind. These are clear limitations, but explicit construction of unitary extension for any commuting tuple of Banach space isometries does not seem to be possible under the shade of orthogonal $2$-sum.

\end{document}